\documentclass[12pt]{article}
\usepackage{latexsym,amsmath,amssymb}
\begin{document}

\baselineskip=20pt

\newcommand{\psp}{\vspace{0.3cm}}

\begin{center}{\LARGE \bf Variational Calculus of Supervariables }
\end{center}
\begin{center}{\LARGE \bf and Related Algebraic Structures}\end{center}

\vspace{0.2cm}

\begin{center}{\large Xiaoping Xu}\end{center}
\begin{center}{Department of Mathematics, The Hong Kong University of Science \& Technology}\end{center}
\begin{center}{Clear Water Bay, Kowloon, Hong Kong}\footnote{Research supported
 by the Direct Allocation Grant 4083 DAG93/94 from HKUST.}\end{center}

\vspace{0.3cm}

\begin{center}{\Large \bf Abstract}\end{center}
\vspace{0.2cm}

{\small We establish a formal variational calculus of supervariables, which is a combination of the bosonic theory of Gel'fand-Dikii and the fermionic theory in our earlier work. Certain interesting new algebraic structures are found in connection with Hamiltonian superoperators in terms of our theory. In particular, we find connections between Hamiltonian superoperators and Novikov-Poisson algebras that we introduced in our earlier work in order to establish a tensor theory of Novikov algebras. Furthermore, we prove that an odd linear Hamiltonian superoperator in our variational calculus induces a Lie superalgebra, which is a natural generalization of the Super-Virasoro algebra under certain conditions.}

\section{Introduction}

Formal variational calculus was introduced by Gel'fand and Dikii [GDi1-2] in studying Hamiltonian systems related to certain nonlinear partial differential equation, such as the KdV equations. Invoking the variational derivatives, they found certain interesting Poisson  structures. Moreover, Gel'fand and Dorfman [GDo] found more connections between Hamiltonian operators and algebraic structures. Balinskii and Novikov [BN] studied similar Poisson structures from another point of view. 

The nature of Gel'fand and Dikii's formal variational calculus is bosonic. In [X3], we presented a general frame of Hamiltonian superoperators and a purely fermionic formal variational calculus. Our work [X3] was based on pure algebraic
analogy. In this paper, we shall present a formal variational calculus of supervariables, which is a combination of the bosonic theory of Gel'fand-Dikii and the fermionic theory in [X3]. Our new theory was motivated by the known super-symmetric theory in mathematical physics (cf. [De], [M]). We find the conditions for a ``matrix differential operator'' to be a Hamiltonian superoperator. In particular, we classify two classes of Hamiltonian superoperators by introducing two kinds of new algebraic structures. Moreover, we prove that an odd linear Hamiltonian superoperator in our variational calculus induces a Lie superalgebra, which is a natural generalization of the Super-Virasoro algebra under certain conditions. We believe that the results in this paper would be useful in study nonlinear super differential equations. They could also play important roles in the application theory of algebras. The discovery of our new algebraic structures proposes new objects in algebraic research. In fact, a new family of infinite-dimensional simple Lie superalgebras were discovered in [X5] based on the results in this paper. 

      Recently, we notice that Daletsky [Da1] introduced a definition of a Hamiltonian superoperator associated with an abstract complex of a Lie superalgebra. He also established in [Da1-2] a formal variational calculus over a commutative superalgebra generated by a set of so-called ``graded  symbols'' with coefficients valued in a Grassman algebra. We believe that one of the subtlenesses of introducing Hamiltonian superoperators is the  constructions of suitable natural complexes of a Lie superalgebra. In our work [X3], we gave a concrete construction of the complex of a colored Lie superalgebra with respect to a graded module and explained the meaning of a Hamiltonian superoperator in detail. It seems to us that the formal variational calculus introduced in [Da1-2] lacks links with the known super-symmetric theory (cf. [De], [M]). For instance, its connection with the known super differential equations, such as the super-symmetric KdV equations, are not clear (cf. [M]). Our formal variational calculus in [X3] was based on free fermionic fields. The combined theory of Gel'fand-Dikki's [GDi1] and ours [X3] that we shall present in this paper  is well motivated by the theory of super-symmetric KdV equations (cf. [M]) and the super-symmetric theory in [De]. Our main purpose in this paper is to show certain new algebraic structures naturally arisen from our theory of Hamiltonian superoperators in a supervariable.

 Below, we shall give more detailed introduction.

Throughout this paper, we let $\Bbb{R}$ be the field of real numbers, and all the vector spaces are assumed over $\Bbb{R}$. Denote by $\Bbb{Z}$ the ring of integers and by $\Bbb{N}$ the set of natural numbers $\{0,1,2,...\}$.
First let us briefly introduce the general frame of Hamiltonian superoperators. We shall sightly modify the differential $d$ defined in (2.7) of [X3]. 

\newcommand{\zt}{\Bbb{Z}_2}

A {\it  Lie superalgebra} $L$ is a $\Bbb{Z}_2$-graded algebra $L=L_0\oplus L_1$ with the operation $[\cdot,\cdot]$ satisfying
$$[x,y]=-(-1)^{xy}[y,x],\qquad [x,[y,z]]+(-1)^{x(y+z)}[y,[z,x]]+(-1)^{z(x+y)}[z,[x,y]]=0\eqno(1.1)$$
for $x,y,z\in L$, where we have used the convention of the notions of exponents of $-1$ used in mathematical physics (cf. [De]); that is, when a vector $u\in L$ appears in an exponent of $-1$, we always means $u\in L_i$ and the value of $u$ in the exponent is $i$. A {\it graded module} $M$ of $L$ is a $\zt$-graded vector space $M=M_0\oplus M_1$ with the action of $L$ on $M$ satisfies:
$$L_i(M_j)\subset M_{i+j},\qquad [x,y]v=xyv-(-1)^{xy}yxv\qquad\mbox{for}\;i,j\in \zt;\;\;x,y\in L;\;v\in M.\eqno(1.2)$$

    A $q$-{\it form of} $L$ {\it with values in} $M$ is a multi-linear map $\omega:\;L^q=L\times \cdots \times L\rightarrow M$ for which
$$\omega (x_1,x_2, \cdots,x_q)=-(-1)^{x_ix_{i+1}}\omega(x_1,\cdots,x_{i-1},x_{i+1},x_i,x_{i+2},\cdots, x_q)\eqno(1.3)$$
for $x_1,...,x_q\in L$. We denote by $c^q(L,M)$ the set of $q$-forms. We define the grading over $c^q(L,M)$ by
$$c^q(L,M)_i=\{\omega \in c^q(L,M)\mid \omega(x_1,...,x_q)\in M_{j_1+\cdots j_q+i}\;\mbox{for}\;x_l\in L_{j_l}\},\qquad i\in \Bbb{Z}_2.\eqno(1.4)$$
Then we have $c^q(L,M)=c^q(L,M)_0+c^q(L,M)_1$.
 Moreover, we define a differential $d:\;c^q(L,M)\rightarrow c^{q+1}(L,M)$ by
\begin{eqnarray*}& & d\omega(x_1,x_2,...,x_{q+1})\\&=&\sum_{i=1}^{q+1}(-1)^{i+1+(\omega+x_1+\cdots x_{i-1})x_i}x_i(\omega(x_1,...,\check{x}_i,...,x_{q+1}))+\sum_{i<j}(-1)^{i+j+(x_1+\cdots +x_{i-1})x_i}\\& &(-1)^{(x_1+\cdots +\check{x}_i+\cdots+ x_{j-1})x_j}\omega([x_i,x_j],x_1,...,\check{x}_i,...,\check{x}_j,...,x_{q+1})\hspace{4cm}(1.5)\end{eqnarray*}
for $\omega\in c^q(L,M),x_l\in L.$ A $q$-form $\omega$ is called {\it closed} if $d\omega=0$. It is easily seen that $d^2=0$ by the proof of Proposition 2.1 in [X3]. 

Let $\omega\in c^2(L,M)_j$. We define:
$${\cal H}_i=\{(x,m)\in L_i\times M_{i+j}\mid \omega (y,x)=(-1)^{jy}ym\;\mbox{for}\;y\in L\},\;\; {\cal H}={\cal H}_0+{\cal H}_1.\eqno(1.6)$$
By (2.10) in [X3], $([x,y],\omega(x,y))\in {\cal H}_{j+l}$ for $x\in L_j,y\in L_l$ if $\omega$ is closed. In this case, we have the following super Poisson bracket
$$\{m_1,m_2\}=\omega(x_1,x_2)\qquad \mbox{for}\;\;(x_1,m_1),(x_2,m_2)\in {\cal H}\eqno(1.7)$$
over the subspace ${\cal N}$ of $M$ defined by
$${\cal N}={\cal N}_0+{\cal N}_1,\qquad {\cal N}_i=\{u\in M_{j+i}\mid (L_i,u)\bigcap {\cal H}\neq\emptyset\}.\eqno(1.8)$$

Let $\Omega$ be a graded subspace of $c^1(L,M)$ such that $dM\subset \Omega$. A graded linear map $H$ is called {\it super skew-symmetric} if
$$\xi_1(H\xi_2)=-(-1)^{(H\xi_1)(H\xi_2)}\xi_2(H\xi_1)\qquad \mbox{for}\;\;\xi_1,\xi_2\in \Omega.\eqno(1.9)$$
With a super skew-symmetric graded linear map $H:\: \Omega\rightarrow L$, we connect a 2-form $\omega_H$ defined on $\mbox{Im}\: H$ by 
$$\omega_H(H\xi_1,H\xi_2)=\xi_2(H\xi_1)\qquad\mbox{for}\;\;\xi_1,\xi_2\in \Omega. \eqno(1.10)$$
We say that $H$ is a {\it Hamiltonian superoperator} if (a) the subspace $\mbox{\it Im}\:H$  of $L$  is a subalgebra; (b) the form $\omega_H$  is closed on $H(\Omega)$.

In [GDo] and [BN],  a new algebra, which was called a ``Novikov algebra'' in [O1], was introduced. A {\it Novikov algebra} ${\cal A}$ is a vector space with an operation ``$\circ$'' satisfying:
$$(x\circ y)\circ z=(x\circ z)\circ y,\qquad (x\circ y)\circ z-x\circ (y\circ z)=(y\circ x)\circ z-y\circ (x\circ z)\eqno(1.11)$$
for $x,y,z\in {\cal A}$. The beauty of a Novikov algebra is that the left multiplication operators forms a Lie algebra and the right multiplication operators are commutative (cf. [Z], [O1]). Zel'manov [Z] proved that any finite-dimensional simple Novikov algebra over an algebraically closed field with characteristic $0$ is one-dimensional. Osborn [O1-5] classified simple Novikov algebras with an idempotent element and their certain modules. In [X4], we gave a complete classification of finite-dimensional simple Novikov algebras and their irreducible modules over an algebraically closed field with prime characteristic. Another algebraic structure introduced in [GDo], which we called ``Gel'fand-Dorfman operator algebra,'' was proved in [X2] to be equivalent to an associative algebra with a derivation under the unitary condition.

 A Novikov algebra actually provides a Poisson structure associated with many-body systems analogous to the KdV-equation (cf. [GDo], [BN]). One might think that the algebra corresponding to the super Poisson structure associated with many-body systems analogous to the super KdV-equations should be the following natural super analogue of Novikov algebras. A {\it Novikov superalgebra} is a $\Bbb{Z}_2$-graded vector space ${\cal A}={\cal A}_0\oplus {\cal A}_1$ with an operation ``$\circ$'' satisfying:
$$(x\circ y)\circ z=(-1)^{yz}(x\circ z)\circ y,\;\; (x\circ y)\circ z-x\circ (y\circ z)=(-1)^{xy}(y\circ x)\circ z-(-1)^{xy}y\circ (x\circ z)\eqno(1.12)$$
for $x,y,z\in {\cal A}$.  It is surprised that Novikov superalgebras are not the algebraic structures corresponding to the super Poisson structures associated with many-body systems analogous to the super KdV-equations. In fact, Novikov superalgebras do not fit in our theory of Hamiltonian superoperators in a supervariable at all. This is because of that the image of a Hamiltonian superoperator is required to be a graded subspace as we introduced in the above. 

As one of the main theorems (see Theorem 3.1), we prove in Section 3 that the algebraic structures corresponding to the Hamiltonian operators (or super Poisson structures) associated with many-body systems (see (3.3)) analogous to the super KdV-equations (see (2.8)) are what we call ``NX-bialgebras.'' An {\it NX-bialgebra} is a vector space $V$ with two operations ``$\times,\circ$'' such that $(V,\times)$ forms a commutative (may not be associative) algebra and $(V,\circ)$ forms a Novikov algebra for which
$$(u\times v)\circ w=u\times (v\circ w),\eqno(1.13)$$
$$(u\times v)\times w+u\times (v\times w)=(v\circ u)\times w+u\times (v\circ w)-v\circ (u\times w),\eqno(1.14)$$
$$(u\times v)\times w-u\times (v\times w)=(u\times v)\circ w+w\circ (u\times v)-u\circ (v\times w)-(v\times w)\circ u\eqno(1.15)$$
for $u,v,w\in V$.

In [X4], we introduced ``Novikov-Poisson'' algebras in order to establish a tensor theory of Novikov algebras. A {\it Novikov-Poisson algebra} is a vector space ${\cal A}$ with two operations ``$\cdot,\circ$'' such that $({\cal A},\cdot)$ forms a commutative associative algebra (may not have an identity element) and $(
{\cal A},\circ)$ forms a Novikov algebra  for which
$$(x\cdot y)\circ z=x\cdot (y\circ z),\qquad (x\circ y)\cdot z-x\circ (y\cdot z)=(y\circ x)\cdot z-y\circ (x\cdot z)\eqno(1.16)$$
for $x,y,z\in {\cal A}$. We prove in Section 3 that certain Novikov-Poisson algebras are NX-bialgebras. This in a way shows the significance of introducing Novikov-Poisson algebras. The detailed study on Novikov Poisson algebras was carried in our work [X5].

We can view the algebraic structure (1.11) as a {\it bosonic Novikov algebra} because the right multiplication operators are commutative. In Section 4, we prove that the following ``fermionic Novikov algebra'' does correspond to a certain Hamiltonian superoperator in a supervariable. 
A {\it fermionic Novikov algebra} ${\cal A}$ is a vector space with an operation ``$\circ$'' satisfying:
$$(x\circ y)\circ z=-(x\circ z)\circ y,\qquad (x\circ y)\circ z-x\circ (y\circ z)=(y\circ x)\circ z-y\circ (x\circ z)\eqno(1.17)$$
for $x,y,z\in {\cal A}$. 

In Section 5, we prove that an odd linear Hamiltonian superoperator induces a Lie superalgebra, which is a natural generalization of the Super-Virasoro algebra under certain conditions. Section 2 is the general theory of our formal variational calculus of supervariables.

\section{Formal Calculus}

In this section, we shall present the frame of our variational calculus of super variables. 

\newcommand{\Lmd}{\Lambda}

Let $\Lambda$ be a vector space that is not necessary finite-dimensional. Let $F(\Lambda)$ be the free associative algebra generated by $\Lmd$. Then the exterior algebra $R$ generated by $\Lmd$ is isomorphic to 
$$ R=F(\Lmd)/(\{uv+vu\mid u,v\in\Lmd\}).\eqno(2.1)$$
We can identify $\Lmd$ with its image in $R$. Note that 
$$R=\Bbb{R}\oplus \Lmd  R=R_c\oplus  R_a,\qquad\mbox{where}\;\;  R_c=\sum_{n=0}^{\infty}\Lmd^{2n},\;\;R_a=\sum_{n=0}^{\infty}\Lmd^{2n+1}.\eqno(2.2)$$
According to [De], the elements of $R_c$ are called $c$-{\it numbers} (means commutative numbers) and  the elements of $R_a$ are called $a$-{\it numbers} (means anti-commutative numbers). Any $u\in R$ can be uniquely written $u=u_b+u_s$ with $u_b\in \Bbb{R},\;u_s\in \Lmd R$ and $u_b$ ($u_s$) is called the {\it body} ({\it soul}, respectively) of $u$. Any analytic function $f$ from $R_c$ to $R$ is of the form
$$f(x)=\sum_{n=0}^{\infty}{\phi^{(n)}(x_b)\over n!}x_s^n,\qquad\mbox{where}\;\;\phi:\Bbb{R}\rightarrow R\;\mbox{is}\;C^{\infty}.\eqno(2.3)$$
An analytic function $\Psi:\;R_c\times R_a \rightarrow R$ is of the form
$$\Psi(x,\theta)=f_0(x)+f_1(x)\theta,\qquad\mbox{where}\;\;f_i:\;R_c\rightarrow R\;\mbox{are analytic}\eqno(2.4)$$
(cf. [De]). Note that
 $$\theta^2=0,\qquad\partial_{\theta}^2=0.\eqno(2.5)$$
 Define 
$$D=\theta\partial_x+\partial_{\theta}\eqno(2.6)$$
Then
$$D^2=\partial_x\eqno(2.7)$$
(cf. [M]). Let $\Phi(x,\theta,t)$ be a function from $R_c\times R_a\times \Bbb{R}$ to $R$. Moreover, we assume that $\Phi(x,\theta, t)\in R_a$ for any $(x,\theta,t)\in (R_c\times R_a\times \Bbb{R})$, $\Phi$ is analytic for fixed $t$ and is $C^1$ with respect to $t$. A super KdV equation is of form
$$\Phi_t=-D^6\Phi+\mu D^2(\Phi D\Phi)+(6-2\mu)D\Phi D^2\Phi\eqno(2.8)$$
(cf. [M]). Mathieu [M] found the Hamiltonians for the above equation when $\mu=2, 3$.

Let $\{\Phi_i\mid I\}$ be a family of functions from $R_c\times R_a\times \Bbb{R}$ to $R$ with the same properties as the above $\Phi$. Set
$$\Phi_i(n+1)=D^n\Phi_i\qquad \mbox{for}\;i\in I;\;n\in \Bbb{N}.\eqno(2.9)$$
Then we have
$$\Phi_i(m)\Phi_j(n)=(-1)^{mn}\Phi_j(n)\Phi_i(m)\qquad \mbox{for}\;i,j\in I;\;m,n\in \Bbb{N}^+=\Bbb{N}\setminus\{0\}.\eqno(2.10)$$
Let ${\cal A}$ be the subalgebra generated by $\{\Phi_i(n)\mid i\in I,\;n\in \Bbb{N}^+\}$ (the set of functions from $R_c\times R_a\times \Bbb{R}$ to $R$ forms an associative algebra). Note that ${\cal A}$ is a $\Bbb{Z}_2$-graded algebra ${\cal A}={\cal A}_0+{\cal A}_1$ with
$${\cal A}_i=\mbox{span}\{\Phi_{i_1}(n_1)\cdots \Phi_{i_p}(n_p)\mid p\in \Bbb{N},i_j\in I,\; n_j\in \Bbb{N}^+, \sum_{j=1}^pn_j\equiv i\;(\mbox{mod}\;2)\},\eqno(2.11)$$
$$u_1u_2=(-1)^{u_1u_2}u_2u_1,\;\;D(u_1u_2)=D(u_1)u_2+(-1)^{u_1}u_1D(u_2)\qquad\mbox{for}\;\;u_1,u_2\in {\cal A}.\eqno(2.12)$$

\newcommand{\ptl}{\partial}

Now we treat $\{\Phi_i(n)\}$ as formal variables. Set
$$L_i=\{\sum_{j\in I}\sum_{l\in \Bbb{N}^+}u_{j,l}\ptl_{\Phi_j(l)}\mid u_{j,l}\in {\cal A}_{i+l}\},\;\;i\in \Bbb{Z}_2,\qquad L=L_0+L_1.\eqno(2.13)$$
Note that the set of superderivations of ${\cal A}$ forms a Lie superalgebra. 
In particular,  $L$ forms a Lie sub-superalgebra with the commutator: 
$$[\partial_1,\partial_2]=\sum_{j,p\in I}\sum_{l,q\in \Bbb{N}^+}(u^1_{p,q}\ptl_{\Phi_p(q)}(u^2_{j,l})-(-1)^{\ptl_1\ptl_2}u^2_{p,q}\ptl_{\Phi_p(q)}(u^1_{j,l}))\partial_{\Phi_j(l)}\eqno(2.14)$$
for $\partial_s=\sum_{j\in I}\sum_{l\in \Bbb{N}^+} u^s_{j,l}\ptl_{\Phi_j(l)}\in L$. 
Note that we can write
$$D=\sum_{i\in I}\sum_{n\in \Bbb{N}^+}\Phi_i(n+1)\ptl_{\Phi_i(n)}.\eqno(2.15)$$
Thus $D\in L$.

By the proof of Lemma 3.2 in [X3], we have:
\psp

{\bf Lemma 2.1}. {\it For} $\partial=\sum_{j\in I}\sum_{l\in \Bbb{N}^+}u_{j,l}\ptl_{\Phi_j(l)}\in (L_0\bigcup L_1)$, $[\partial, D]=0$ {\it if and only if}
$$u_{j,n+1}=(-1)^{n\ptl}D^n(u_{j,1}),\qquad  n\in \Bbb{N}.\eqno(2.16)$$
\psp

Set
$${\cal L}={\cal L}_1+{\cal L}_2\subset {\cal A}^I,\qquad {\cal L}_s=({\cal A}_{s+1})^I.\eqno(2.17)$$
For any $\bar{u}=\{u_i\mid i\in I\}\in {\cal L}_s$, we let
$$\partial_{\bar{u}}=\sum_{j\in I}\sum_{n\in \Bbb{N}}(-1)^{sn}D^n(u_j)\ptl_{\Phi_j(n+1)}\in L.\eqno(2.18)$$
Then $[\partial_{\bar{u}},D]=0$. 

For $\bar{u}=\{u_i\},\bar{v}=\{v_i\}\in {\cal L}$,
$$[\partial_{\bar{u}},\partial_{\bar{v}}]=\ptl_{\bar{w}}\eqno(2.19)$$
with

\begin{eqnarray*}\hspace{1cm}\bar{w}&=&\{\sum_{p\in I}\sum_{ m\in \Bbb{N}^+}((-1)^{m\bar{u}}D^m(u_p)\ptl_{\Phi_p(m+1)}(v_q)\\& &-(-1)^{\bar{u}\bar{v}+m\bar{v}}D^m(v_p)\ptl_{\Phi_p(m+1)}(u_q))\mid q\in I\}\hspace{5.2cm}(2.20)\end{eqnarray*}
(cf. (3.26-27) in [X3]). 

Thus if we  define
$$[\bar{u},\bar{v}]=\bar{w},\eqno(2.21)$$
then $({\cal L},\Bbb{Z}_2,[\cdot,\cdot])$ forms a Lie superalgebra. 

Next we define variational operators on ${\cal A}$:
$$\delta_i=\sum_{m=0}^{\infty}(-1)^{m(m-1)/2}D^m\circ \ptl_{\Phi_i(m+1)},\qquad \bar{\delta}=\{\delta_i\mid i\in I\}.\eqno(2.22)$$

By the proof of Lemma 3.4 in [X3], we have:
\psp

{\bf Lemma 2.2}. {\it For any} $u\in \sum_{i\in I,n\in \Bbb{N}}{\cal A}\Phi_i(n+1)$, 
$$\bar{\delta}(u)=0\Longleftrightarrow u=D(v)\;\;\mbox{\it for some}\;\;v\in {\cal A}.\eqno(2.23)$$
\psp

Now we let
$$\tilde{\cal A}={\cal A}/D({\cal A}).\eqno(2.24)$$
We define an action of ${\cal L}$ on $\tilde{\cal A}$ by
$$\bar{u}(\tilde{w})=\partial_{\bar{u}}(w)+D({\cal A}
)= \sum_{i\in I}(u_i\delta_i(w))^{\sim}\eqno(2.25)$$
(cf. (3.39) in [X3]). 
This is well defined since $[\partial_{\bar{u}},D]=0$. Thus $\tilde{\cal A}$ forms an ${\cal L}$-module. Furthermore, we set 
$$\Omega=\{\bar{\xi}=\{\xi_i\}\in {\cal A}^I\mid \mbox{only finite number of}\;\xi_i\neq 0\}.\eqno(2.26)$$
For any $\bar{\xi}\in \Omega, \;\bar{u}\in {\cal L}$, we define:
$$\bar{\xi}(\bar{u})=\sum_{i\in I}(u_i\xi_i)^{\sim}.\eqno(2.27)$$
Then $\Omega\subset c^1({\cal L},\tilde{\cal A})$. Note that by (2.25), 
$$d(\tilde{w})=\bar{\delta}(w)\in \Omega\qquad\mbox{for}\;\;\tilde{w}\in {\cal A},\eqno(2.28)$$
where (2.23) implies that the map $\bar{\delta}:\:\tilde{\cal A}\rightarrow \Omega$ is well defined. Hence $d(\tilde{\cal A})\in \Omega$.

Note that as sets, $\Omega\subset {\cal L}$. We let
$$\Omega_i=\Omega\bigcap {\cal L}_i\qquad\mbox{for}\;\;i\in \Bbb{Z}_2.\eqno(2.29)$$
Suppose that $H:\;\Omega\rightarrow {\cal L}$ is a linear map as follows: for $\bar{\xi}\in \Omega_i,\;i\in \Bbb{Z}_2$,
$$(H\bar{\xi})_p=\sum_{q\in I}H^i_{p,q}\xi_q,\;\;\mbox{where}\;\;H^i_{p,q}=\sum_{l=0}^{n(i,p,q)}a_{p,q,l}^iD^l\;\;\mbox{with}\;\;a_{p,q,l}^i\in{\cal A}_{\iota+l},\;\iota\in \Bbb{Z}_2.\eqno(2.30)$$
Such an $H$ is called a {\it matrix differential operator of type} $\iota$. Moreover, $H(\Omega)$ is a $\Bbb{Z}_2$-graded subspace. Furthermore, the super skew-symmetry is equivalent to
$$\sum_{l=0}^{n(0,p,q)}(-1)^{(2\iota+l)(l-1)/2}D^l\circ a_{p,q,l}^0=\sum_{l=0}^{n(0,q,p)}a_{q,p,l}^0D^l,\;\;\;a_{p,q,l}^0=(-1)^{\iota+1}a_{p,q,l}^1.\eqno(2.31)$$

Let $H:\;\Omega \rightarrow {\cal L}$ be a super skew-symmetric matrix differential operator. We want to find  the condition for $H$ to be a Hamiltonian operator. For $\bar{\xi}\in \Omega_i$, we define a linear map $(D_H\bar{\xi}):\;{\cal L}\rightarrow {\cal L}$ by
$$(D_H\bar{\xi})(\bar{\eta})=(D_H\bar{\xi})\bar{\eta},\;\;\;(D_H\bar{\xi})_{p,q}=\sum_{t\in I}\sum_{l,m\in \Bbb{N}}(-1)^{m(i+\iota)}\ptl_{\Phi_q(m+1)}(a^i_{p,t,l})D^l(\xi_t)D^m, \eqno(2.32)$$
for $\bar{\eta}\in \Omega$. 

By the proof of Theorem 4.1 in [X3], we have:
\psp

{\bf Theorem 2.3}. {\it A matrix differential operator} $H$ {\it of form (2.30) is a Hamiltonian operator if and only if (2.31) and the following equation hold:}
\begin{eqnarray*} & &(-1)^{\bar{\xi}_1}\bar{\xi}_3((D_H\bar{\xi}_1)H\bar{\xi}_2)+(-1)^{\bar{\xi}_2+(\bar{\xi}_1+\iota,\bar{\xi}_2+\bar{\xi}_3)}\bar{\xi}_1((D_H\bar{\xi}_2)H\bar{\xi}_3)\\&=&-(-1)^{\bar{\xi}_3+(\bar{\xi}_3+\iota,\bar{\xi}_1+\bar{\xi}_2)}\bar{\xi}_2((D_H\bar{\xi}_3)H\bar{\xi}_1)\hspace{7.8cm}(2.33)\end{eqnarray*}
{\it for} $\bar{\xi}_1,\bar{\xi}_2,\bar{\xi}_3\in \Omega$.
\psp

{\bf Remark 2.4}. By (2.31) and the above theorem, the operator
$$H=\sum_{m=0}^{\infty} a_mD^{4m+1},\;\;\;a_m\in \Bbb{R}\eqno(2.34)$$
is a Hamiltonian operator of type 1. Moreover, the operator $H'$ defined by
$$H'(\bar{\xi})=(-1)^{\bar{\xi}}\sum_{m=0}^{\infty} b_mD^{4m}\bar{\xi},\;\;\;b_m\in \Bbb{R},\;\;\mbox{for}\;\;\bar{\xi}\in \Omega,\eqno(2.35)$$
is a Hamiltonian operator of type 0.
\psp

Let $H_1$ and $H_2$ be matrix differential operators of the same type $\iota$. If $aH_1+bH_2$ is Hamiltonian for any $a,b\in \Bbb{R}$, then we call $(H_1,H_2)$ a {\it Hamiltonian pair}. 
For any two matrix differential operators $H_1$ and $H_2$, we define the Schouten-Nijenhuis super-bracket $[H_1,H_2]:\: \Omega^3\rightarrow \tilde{\cal A}$ by
\begin{eqnarray*}& &[H_1,H_2](\bar{\xi}_1,\bar{\xi}_2,\bar{\xi}_3)\\&=&
(-1)^{\bar{\xi}_1}\bar{\xi}_3((D_{H_1}\bar{\xi}_1)H_2\bar{\xi}_2)+(-1)^{\bar{\xi}_1}\bar{\xi}_3((D_{H_2}\bar{\xi}_1)H_1\bar{\xi}_2)\\& &+(-1)^{\bar{\xi}_2+(\bar{\xi}_1+\iota,\bar{\xi}_2+\bar{\xi}_3)}\bar{\xi}_1((D_{H_1}\bar{\xi}_2)H_2\bar{\xi}_3)+(-1)^{\bar{\xi}_2+(\bar{\xi}_1+\iota,\bar{\xi}_2+\bar{\xi}_3)}\bar{\xi}_1((D_{H_2}\bar{\xi}_2)H_1\bar{\xi}_3)\\& &+(-1)^{\bar{\xi}_3+(\bar{\xi}_3+\iota,\bar{\xi}_1+\bar{\xi}_2)}\bar{\xi}_2((D_{H_1}\bar{\xi}_3)H_2\bar{\xi}_1)+(-1)^{\bar{\xi}_3+(\bar{\xi}_3+\iota,\bar{\xi}_1+\bar{\xi}_2)}\bar{\xi}_2((D_{H_2}\bar{\xi}_3)H_1\bar{\xi}_1)\hspace{1cm}(2.36)\end{eqnarray*}
for $\bar{\xi}_1,\bar{\xi}_2,\bar{\xi}_3\in \Omega$. Then (2.33) is equivalent to $[H,H]=0$. In general, we have:
\psp

{\bf Corollary 2.5}. {\it  Matrix differential operators} $H_1$ {\it and} $H_2$ {\it of the same type forms a Hamiltonian pair if and only if  they satisfy (2.31) and}
$$[H_1,H_1]=0,\;\;\;[H_2,H_2]=0,\;\;\;[H_1,H_2]=0.\eqno(2.37)$$

\section{Hamiltonian Superoperators and NX-Bialgebras}

\newcommand{\al}{\alpha}
\newcommand{\be}{\beta}
\newcommand{\gm}{\gamma}
\newcommand{\la}{\langle}
\newcommand{\ra}{\rangle}
\newcommand{\lmd}{\lambda}

In this section, we consider the type-1 Hamiltonian operator $H$ of the form:
$$H^1_{\al,\be}=H^0_{\al,\be}=a_{\al,\be}D^5+\sum_{\gm\in I}[b^{\gm}_{\al,\be}\Phi_{\gm}D^2+c_{\al,\be}^{\gm}\Phi_{\gm}(2)D+d_{\al,\be}^{\gm}\Phi_{\gm}(3)],\eqno(3.1)$$
where $a_{\al,\be}^{\gm},b_{\al,\be}^{\gm},c_{\al,\be}^{\gm},d_{\al,\be}^{\gm}\in \Bbb{R}.$ We let 
$$L=\sum_{\al,\be\in I}\chi_{\al,\be}\Phi_{\al}\Phi_{\be}(2),\qquad \chi_{\al,\be}\in \Bbb{R}.\eqno(3.2)$$
Then we have the following many-body systems analogous to  the super KdV equations:
$$(\Phi_{\al})_t=\sum_{\be\in I} H_{\al,\be}\delta_{\be}(L),\qquad\qquad\al\in I.\eqno(3.3)$$

As we shall show below, it is not easy to find the condition for an operator in (3.1) to be Hamiltonian. The difficulty is that (2.33) is equivalent to a set of many equations. Therefore, high technical reductions are needed in order to find the condition of simplest form.

Note that the super skew-symmetry of $H$ is equivalent to
\begin{eqnarray*}& &a_{\al,\be}D^5+\sum_{\gm\in I}[b^{\gm}_{\al,\be}\Phi_{\gm}D^2+c_{\al,\be}^{\gm}\Phi_{\gm}(2)D+d_{\al,\be}^{\gm}\Phi_{\gm}(3)]\\&=&a_{\be,\al}D^5+\sum_{\gm\in I}[b^{\gm}_{\be,\al}D^2\circ\Phi_{\gm}+c_{\be,\al}^{\gm}D\circ\Phi_{\gm}(2)-d_{\be,\al}^{\gm}\Phi_{\gm}(3)]\\&=& a_{\be,\al}D^5+\sum_{\gm\in I}[b^{\gm}_{\be,\al}\Phi_{\gm}D^2+c_{\be,\al}^{\gm}\Phi_{\gm}(2)D+(b_{\be,\al}^{\gm}+c_{\be,\al}^{\gm}-d_{\be,\al}^{\gm})\Phi_{\gm}(3)]\hspace{2.6cm}(3.4)\end{eqnarray*}
by(2.31), equivalently,
$$a_{\al,\be}=a_{\be,\al},\;\;b_{\al,\be}^{\gm}=b_{\be,\al}^{\gm},\;\;c_{\al,\be}^{\gm}=c_{\be,\al}^{\gm},\;\;b_{\al,\be}^{\gm}+c_{\al,\be}^{\gm}=d^{\gm}_{\al,\be}+d^{\gm}_{\be,\al}.\eqno(3.5)$$
Moreover, we let 
$$V=\sum_{\al\in I}\Bbb{R}\Phi_{\al}\eqno(3.6)$$
and define the operations: $\cdot,\times,\circ:\;V\times V\rightarrow V$ and the bilinear form $\la\cdot,\cdot\ra$ by
$$\Phi_{\al}\cdot \Phi_{\be}=\sum_{\gm\in I}b_{\al,\be}^{\gm}\Phi_{\gm},\;\;\Phi_{\al}\times \Phi_{\be}=\sum_{\gm\in I}c_{\al,\be}^{\gm}\Phi_{\gm},\;\;\Phi_{\al}\circ \Phi_{\be}=\sum_{\gm\in I}d_{\al,\be}^{\gm}\Phi_{\gm},\;\;\la \Phi_{\al}, \Phi_{\be}\ra=a_{\al,\be}\eqno(3.7)$$
for $\al,\;\be\in I.$ Then $(V,\cdot),\;(V,\times)$ are commutative algebras (may not be associative) and $\la\cdot,\cdot\ra$ is a symmetric bilinear form. 

In order to find the conditons for which (2.33) holds, we have to find the exact formula for each term in (3.33). For $\bar{\xi}_1,\bar{\xi}_2,\bar{\xi}_3\in \Omega$,  we have
\begin{eqnarray*}& &\bar{\xi}_3((D_H\bar{\xi}_1)H\bar{\xi}_2)\\
&=& \sum_{\al,\be,\gm,\lmd,\mu\in I}[b_{\gm,\al}^{\lmd}D^2(\xi_{1\al})+(-1)^{\bar{\xi}_1+1}c_{\gm,\al}^{\lmd}D(\xi_{1\al})D+d_{\gm,\al}^{\lmd}\xi_{1\al}D^2]\\& &
[a_{\lmd,\be}D^5(\xi_{2\be})+b_{\lmd,\be}^{\mu}\Phi_{\mu}D^2(\xi_{2\be})+c_{\lmd,\be}^{\mu}\Phi_{\mu}(2)D(\xi_{2\be})+d_{\lmd,\be}^{\mu}\Phi_{\mu}(3)\xi_{2\be}]\xi_{3\gm}\\&=& \sum_{\al,\be,\gm\in I}\{\la \Phi_{\gm}\cdot \Phi_{\al}, \Phi_{\be}\ra D^2(\xi_{1\al})D^5(\xi_{2\be})+(-1)^{\bar{\xi}_1+1}[(\Phi_{\gm}\cdot \Phi_{\al})\cdot \Phi_{\be}]D^2(\xi_{1\al})D^2(\xi_{2\be})\\& &+[(\Phi_{\gm}\cdot \Phi_{\al})\times \Phi_{\be}](2)D^2(\xi_{1\al})D(\xi_{2\be})+(-1)^{\bar{\xi}_1+1}[(\Phi_{\gm}\cdot \Phi_{\al})\circ \Phi_{\be}](3)D^2(\xi_{1\al})\xi_{2\be}\\& &+(-1)^{\bar{\xi}_1+1}\la\Phi_{\gm}\times\Phi_{\al}, \Phi_{\be}\ra D(\xi_{1\al})D^6(\xi_{2\be})+(-1)^{\bar{\xi}_1+1}[(\Phi_{\gm}\times \Phi_{\al})\cdot \Phi_{\be}](2)D(\xi_{1\al})D^2(\xi_{2\be})\\& &+[(\Phi_{\gm}\times \Phi_{\al})\cdot \Phi_{\be}]D(\xi_{1\al})D^3(\xi_{2\be})+(-1)^{\bar{\xi}_1+1}[(\Phi_{\gm}\times \Phi_{\al})\times \Phi_{\be}](2)D(\xi_{1\al})D^2(\xi_{2\be})\hspace{3cm}\end{eqnarray*}
\begin{eqnarray*}&& -[(\Phi_{\gm}\times \Phi_{\al})\times \Phi_{\be}](3)D(\xi_{1\al})D(\xi_{2\be})+(-1)^{\bar{\xi}_1+1}[(\Phi_{\gm}\times \Phi_{\al})\circ \Phi_{\be}](4)D(\xi_{1\al})\xi_{2\be}\\& &+[(\Phi_{\gm}\times \Phi_{\al})\circ \Phi_{\be}](3)D(\xi_{1\al})D(\xi_{2\be})+\la \Phi_{\gm}\circ \Phi_{\al}, \Phi_{\be}\ra\xi_{1\al}D^7(\xi_{2\be})\\&&+(-1)^{\bar{\xi}_1+1}[(\Phi_{\gm}\circ \Phi_{\al})\cdot \Phi_{\be}](3)\xi_{1\al}D^2(\xi_{2\be})+(-1)^{\bar{\xi}_1+1}[(\Phi_{\gm}\circ \Phi_{\al})\cdot \Phi_{\be}]\xi_{1\al}D^4(\xi_{2\be})\\& &+[(\Phi_{\gm}\circ \Phi_{\al})\times \Phi_{\be}](4)\xi_{1\al}D(\xi_{2\be})+[(\Phi_{\gm}\circ \Phi_{\al})\times \Phi_{\be}](2)\xi_{1\al}D^3(\xi_{2\be})\\& &+(-1)^{\bar{\xi}_1+1}[(\Phi_{\gm}\circ\Phi_{\al})\circ \Phi_{\be}](5)\xi_{1\al}\xi_{2\be}+(-1)^{\bar{\xi}_1+1}[(\Phi_{\gm}\circ \Phi_{\al})\circ \Phi_{\be}](3)\xi_{1\al}D^2(\xi_{2\be})\}\xi_{3\gm},\hspace{0.4cm}(3.8)\end{eqnarray*}
\begin{eqnarray*}& &(-1)^{(\bar{\xi}_1+1)(\bar{\xi}_2+\bar{\xi}_3)}\bar{\xi}_1((D_H\bar{\xi}_2)H\bar{\xi}_3)\\&=&\sum_{\al,\be,\gm\in I}\{(-1)^{\bar{\xi}_2}\la \Phi_{\al}\cdot \Phi_{\be}, \Phi_{\gm}\ra[D^5(\xi_{1\al})D^2(\xi_{2\be})+2D^3(\xi_{1\al})D^4(\xi_{2\be})+D^1(\xi_{1\al})D^6(\xi_{2\be})\\& &+(-1)^{\bar{\xi}_1+1}D^4(\xi_{1\al})D^3(\xi_{2\be})+(-1)^{\bar{\xi}_1+1}2D^2(\xi_{1\al})D^5(\xi_{2\be})+(-1)^{\bar{\xi}_1+1}\xi_{1\al}D^7(\xi_{2\be})]\\& &+(-1)^{\bar{\xi}_2}[[\Phi_{\al}\cdot \Phi_{\be})\cdot \Phi_{\gm}](3)\xi_{1\al}D^2(\xi_{2\be})+[\Phi_{\al}\cdot \Phi_{\be})\cdot \Phi_{\gm}]D^2(\xi_{1\al})D^2(\xi_{2\be})\\& &+[\Phi_{\al}\cdot \Phi_{\be})\cdot \Phi_{\gm}]\xi_{1\al}D^4(\xi_{2\be})]]+(-1)^{\bar{\xi}_2}[[(\Phi_{\al}\cdot \Phi_{\be})\times \Phi_{\gm}](3)\xi_{1\al}D^2(\xi_{2\be})\\&& +[(\Phi_{\al}\cdot \Phi_{\be})\times \Phi_{\gm}](2)D(\xi_{1\al})D^2(\xi_{2\be})+(-1)^{\bar{\xi}_1+1}[(\Phi_{\al}\cdot \Phi_{\be})\times \Phi_{\gm}](2)\xi_{1\al}D^3(\xi_{2\be})]\\& &
+(-1)^{\bar{\xi}_2+1}[(\Phi_{\al}\cdot \Phi_{\be})\circ \Phi_{\gm}](3)\xi_{1\al}D^2(\xi_{2\be})+(-1)^{\bar{\xi}_1+\bar{\xi}_2+1}\la\Phi_{\al}\times\Phi_{\be}, \Phi_{\gm}\ra[D^6(\xi_{1\al})D(\xi_{2\be})\\& &+3D^4(\xi_{1\al})D^3(\xi_{2\be})+3D^2(\xi_{1\al})D^5(\xi_{2\be})+\xi_{1\al}D^7(\xi_{2\be})]\\ & &
+(-1)^{\bar{\xi}_1+\bar{\xi}_2+1}[[(\Phi_{\al}\times \Phi_{\be})\cdot \Phi_{\gm}](4)\xi_{1\al}D(\xi_{2\be})+[(\Phi_{\al}\times \Phi_{\be})\cdot \Phi_{\gm}](2)D^2(\xi_{1\al})D(\xi_{2\be})\\& &+[(\Phi_{\al}\times \Phi_{\be})\cdot \Phi_{\gm}](2)\xi_{1\al}D^3(\xi_{2\be})]+(-1)^{\bar{\xi}_1+\bar{\xi}_2}[[(\Phi_{\al}\times \Phi_{\be})\cdot \Phi_{\gm}](4)\xi_{1\al}D(\xi_{2\be})\\& &+[(\Phi_{\al}\times \Phi_{\be})\cdot \Phi_{\gm}](2)D^2(\xi_{1\al})D(\xi_{2\be})+[(\Phi_{\al}\times \Phi_{\be})\cdot \Phi_{\gm}](2)\xi_{1\al}D^3(\xi_{2\be})
\\&&-[(\Phi_{\al}\times \Phi_{\be})\cdot \Phi_{\gm}](3)D(\xi_{1\al})D(\xi_{2\be})+(-1)^{\bar{\xi}_1}[(\Phi_{\al}\times \Phi_{\be})\cdot \Phi_{\gm}](3)\xi_{1\al}D^2(\xi_{2\be})\\& &-[(\Phi_{\al}\times \Phi_{\be})\cdot \Phi_{\gm}]D^3(\xi_{1\al})D(\xi_{2\be})-[(\Phi_{\al}\times \Phi_{\be})\cdot \Phi_{\gm}]D(\xi_{1\al})D^3(\xi_{2\be})\\& &+(-1)^{\bar{\xi}_1}[(\Phi_{\al}\times \Phi_{\be})\cdot \Phi_{\gm}]D^2(\xi_{1\al})D^2(\xi_{2\be})+(-1)^{\bar{\xi}_1}[(\Phi_{\al}\times \Phi_{\be})\cdot \Phi_{\gm}]\xi_{1\al}D^4(\xi_{2\be})]\\& &+(-1)^{\bar{\xi}_1+\bar{\xi}_2+1}[[(\Phi_{\al}\times \Phi_{\be})\times \Phi_{\gm}](4)\xi_{1\al}D(\xi_{2\be})+[(\Phi_{\al}\times \Phi_{\be})\times \Phi_{\gm}](2)D^2(\xi_{1\al})D(\xi_{2\be})\\& &+[(\Phi_{\al}\times \Phi_{\be})\times \Phi_{\gm}](2)\xi_{1\al}D^3(\xi_{2\be})]+(-1)^{\bar{\xi}_1+\bar{\xi}_2}[[(\Phi_{\al}\times \Phi_{\be})\times \Phi_{\gm}](4)\xi_{1\al}D(\xi_{2\be})
\\& &-[(\Phi_{\al}\times \Phi_{\be})\times \Phi_{\gm}](3)D(\xi_{1\al})D(\xi_{2\be})+(-1)^{\bar{\xi}_1}[(\Phi_{\al}\times \Phi_{\be})\times \Phi_{\gm}](3)\xi_{1\al}D^2(\xi_{2\be})]\\& &+(-1)^{\bar{\xi}_1+\bar{\xi}_2}[(\Phi_{\al}\times \Phi_{\be})\circ \Phi_{\gm}](4)\xi_{1\al}D(\xi_{2\be})+(-1)^{\bar{\xi}_1+\bar{\xi}_2+1}[[(\Phi_{\al}\times \Phi_{\be})\circ \Phi_{\gm}](4)\xi_{1\al}D(\xi_{2\be})\\&& -[(\Phi_{\al}\times \Phi_{\be})\circ \Phi_{\gm}](3)D(\xi_{1\al})D(\xi_{2\be})+(-1)^{\bar{\xi}_1}[(\Phi_{\al}\times \Phi_{\be})\circ \Phi_{\gm}](3)\xi_{1\al}D^2(\xi_{2\be})]\\& &+(-1)^{\bar{\xi}_2+1}\la \Phi_{\al}\circ \Phi_{\be}, \Phi_{\gm}\ra[D^7(\xi_{1\al})\xi_{2\be}+3D^5(\xi_{1\al})D^2(\xi_{2\be})+3D^3(\xi_{1\al})D^4(\xi_{2\be})\\& &+D(\xi_{1\al})D^6(\xi_{2\be})+(-1)^{\bar{\xi}_1+1}[D^6(\xi_{1\al})D(\xi_{2\be})+3D^4(\xi_{1\al})D^3(\xi_{2\be})\\& &+3D^2(\xi_{1\al})D^5(\xi_{2\be})+\xi_{1\al}D^7(\xi_{2\be})]]+(-1)^{\bar{\xi}_2}[[(\Phi_{\al}\circ \Phi_{\be})\cdot \Phi_{\gm}](5)\xi_{1\al}\xi_{2\be}\\& &+[(\Phi_{\al}\circ \Phi_{\be})\cdot \Phi_{\gm}](3)D^2(\xi_{1\al})\xi_{2\be}+[(\Phi_{\al}\circ \Phi_{\be})\cdot \Phi_{\gm}](3)\xi_{1\al}D^2(\xi_{2\be})]\\& &+(-1)^{\bar{\xi}_2+1}[[(\Phi_{\al}\circ \Phi_{\be})\cdot \Phi_{\gm}](5)\xi_{1\al}\xi_{2\be}+2[(\Phi_{\al}\circ \Phi_{\be})\cdot \Phi_{\gm}](3)D^2(\xi_{1\al})\xi_{2\be}\hspace{5cm}\end{eqnarray*}
\begin{eqnarray*}& &+2[(\Phi_{\al}\circ \Phi_{\be})\cdot \Phi_{\gm}](3)\xi_{1\al}D^2(\xi_{2\be})+2[(\Phi_{\al}\circ \Phi_{\be})\cdot \Phi_{\gm}]D^2(\xi_{1\al})D^2(\xi_{2\be})\\& &+[(\Phi_{\al}\circ \Phi_{\be})\cdot \Phi_{\gm}]D^4(\xi_{1\al})\xi_{2\be}+[(\Phi_{\al}\circ \Phi_{\be})\cdot \Phi_{\gm}]\xi_{1\al}D^4(\xi_{2\be})]
\\& & +(-1)^{\bar{\xi}_2}[[(\Phi_{\al}\circ \Phi_{\be})\times \Phi_{\gm}](5)\xi_{1\al}\xi_{2\be}+[(\Phi_{\al}\circ \Phi_{\be})\times \Phi_{\gm}](4)D(\xi_{1\al})\xi_{2\be}\\& &+(-1)^{\bar{\xi}_1+1}[(\Phi_{\al}\circ \Phi_{\be})\times \Phi_{\gm}](4)\xi_{1\al}D(\xi_{2\be})]+(-1)^{\bar{\xi}_2+1}[[(\Phi_{\al}\circ \Phi_{\be})\times \Phi_{\al}](5)\xi_{1\al}\xi_{2\be}\\& &+[(\Phi_{\al}\circ \Phi_{\be})\times \Phi_{\al}](3)D^2(\xi_{1\al})\xi_{2\be}+[(\Phi_{\al}\circ \Phi_{\be})\times \Phi_{\al}](3)\xi_{1\al}D^2(\xi_{2\be})\\& &+[(\Phi_{\al}\circ \Phi_{\be})\times \Phi_{\al}](4)D(\xi_{1\al})\xi_{2\be}+(-1)^{\bar{\xi}_1+1}[(\Phi_{\al}\circ \Phi_{\be})\times \Phi_{\al}](4)\xi_{1\al}D(\xi_{2\be})\\& &+[(\Phi_{\al}\circ \Phi_{\be})\times \Phi_{\al}](2)D^3(\xi_{1\al})\xi_{2\be}+[(\Phi_{\al}\circ \Phi_{\be})\times \Phi_{\al}](2)D(\xi_{1\al})D^2(\xi_{2\be})
\\& &+(-1)^{\bar{\xi}_1+1}[(\Phi_{\al}\circ \Phi_{\be})\times \Phi_{\al}](2)D^2(\xi_{1\al})D(\xi_{2\be})+(-1)^{\bar{\xi}_1+1}[(\Phi_{\al}\circ \Phi_{\be})\times \Phi_{\al}](2)\xi_{1\al}D^3(\xi_{2\be})]\\& &+(-1)^{\bar{\xi}_2+1}[(\Phi_{\al}\circ\Phi_{\be})\circ \Phi_{\gm}](5)\xi_{1\al}\xi_{2\be}+(-1)^{\bar{\xi}_2}[[(\Phi_{\al}\circ \Phi_{\be})\circ \Phi_{\gm}](5)\xi_{1\al}\xi_{2\be}\\& &+[(\Phi_{\al}\circ \Phi_{\be})\circ \Phi_{\gm}](3)D^2(\xi_{1\al})\xi_{2\be}+[(\Phi_{\al}\circ \Phi_{\be})\circ \Phi_{\gm}](3)\xi_{1\al}D^2(\xi_{2\be})]]\}\xi_{3\gm},\hspace{1.7cm}(3.9)\end{eqnarray*}
\begin{eqnarray*}& &(-1)^{(\bar{\xi}_3+1)(\bar{\xi}_1+\bar{\xi}_2)}\bar{\xi}_2((D_H\bar{\xi}_3)H\bar{\xi}_1)\\&=&
\sum_{\al,\be,\gm\in I}\{(-1)^{\bar{\xi}_3}\la \Phi_{\be}\cdot \Phi_{\gm}, \Phi_{\al}\ra[D^7(\xi_{1\al}) \xi_{2\be}+D^5(\xi_{1\al})D^2(\xi_{2\be})]
\\& &+(-1)^{\bar{\xi}_3}[[(\Phi_{\be}\cdot \Phi_{\gm})\cdot \Phi_{\al}](3)D^2(\xi_{1\al})\xi_{2\be}+[(\Phi_{\be}\cdot \Phi_{\gm})\cdot \Phi_{\al}]D^4(\xi_{1\al})\xi_{2\be}\\& &+[(\Phi_{\be}\cdot \Phi_{\gm})\cdot \Phi_{\al}]D^2(\xi_{1\al})D^2(\xi_{2\be})]+(-1)^{\bar{\xi}_3}[[(\Phi_{\be}\cdot \Phi_{\gm})\times \Phi_{\al}](4)D(\xi_{1\al})\xi_{2\be}\\& &+[(\Phi_{\be}\cdot \Phi_{\gm})\times \Phi_{\al}](2)D^3(\xi_{1\al})\xi_{2\be}+[(\Phi_{\be}\cdot \Phi_{\gm})\times \Phi_{\al}](2)D(\xi_{1\al})D^2(\xi_{2\be})]\\& &+(-1)^{\bar{\xi}_3}[[(\Phi_{\be}\cdot \Phi_{\gm})\circ \Phi_{\al}](5)\xi_{1\al}\xi_{2\be}+[(\Phi_{\be}\cdot \Phi_{\gm})\circ \Phi_{\al}](3)D^2(\xi_{1\al})\xi_{2\be}\\& &+[(\Phi_{\be}\cdot \Phi_{\gm})\circ \Phi_{\al}](3)\xi_{1\al}D^2(\xi_{2\be})]+(-1)^{\bar{\xi}_3}\la\Phi_{\be}\times\Phi_{\gm}, \Phi_{\al}\ra [D^7(\xi_{1\al})\xi_{2\be}\\& &+(-1)^{\bar{\xi}_1+1}D^6(\xi_{1\al})D(\xi_{2\be})]+(-1)^{\bar{\xi}_3}[[(\Phi_{\be}\times \Phi_{\gm})\cdot \Phi_{\al}](3)D^2(\xi_{1\al})\xi_{2\be}\\& &+[(\Phi_{\be}\times \Phi_{\gm})\cdot \Phi_{\al}](2)D^3(\xi_{1\al})\xi_{2\be}+(-1)^{\bar{\xi}_1+1}[(\Phi_{\be}\times \Phi_{\gm})\cdot \Phi_{\al}](2)D^2(\xi_{1\al})D(\xi_{2\be})]
\\& &+(-1)^{\bar{\xi}_3+1}[[(\Phi_{\be}\times \Phi_{\gm})\cdot \Phi_{\al}](2)D^3(\xi_{1\al})\xi_{2\be}-[(\Phi_{\be}\times \Phi_{\gm})\cdot \Phi_{\al}]D^4(\xi_{1\al})\xi_{2\be}\\& &+(-1)^{\bar{\xi}_1+1}[(\Phi_{\be}\times \Phi_{\gm})\cdot \Phi_{\al}]D^3(\xi_{1\al})D(\xi_{2\be})]
+(-1)^{\bar{\xi}_3}[[(\Phi_{\be}\times \Phi_{\gm})\times \Phi_{\al}](3)D^2(\xi_{1\al})\xi_{2\be}\\& &+[(\Phi_{\be}\times \Phi_{\gm})\times \Phi_{\al}](2)D^3(\xi_{1\al})\xi_{2\be}+(-1)^{\bar{\xi}_1+1}[(\Phi_{\be}\times \Phi_{\gm})\times \Phi_{\al}](2)D^2(\xi_{1\al})D(\xi_{2\be})]
\\&&+(-1)^{\bar{\xi}_3}[[(\Phi_{\be}\times \Phi_{\gm})\times \Phi_{\al}](4)D(\xi_{1\al})\xi_{2\be}-[(\Phi_{\be}\times \Phi_{\gm})\times \Phi_{\al}](3)D^2(\xi_{1\al})\xi_{2\be}\\& &+(-1)^{\bar{\xi}_1+1}[(\Phi_{\be}\times \Phi_{\gm})\times \Phi_{\al}](3)D(\xi_{1\al})D(\xi_{2\be})]+(-1)^{\bar{\xi}_3}[[(\Phi_{\be}\times \Phi_{\gm})\circ \Phi_{\al}](5)\xi_{1\al}\xi_{2\be}\\& &+[(\Phi_{\be}\times \Phi_{\gm})\circ \Phi_{\al}](4)D(\xi_{1\al})\xi_{2\be}+(-1)^{\bar{\xi}_1+1}[(\Phi_{\be}\times \Phi_{\gm})\circ \Phi_{\al}](4)\xi_{1\al}D(\xi_{2\be})]
\\& &+(-1)^{\bar{\xi}_3+1}[[(\Phi_{\be}\times \Phi_{\gm})\circ \Phi_{\al}](4)D(\xi_{1\al})\xi_{2\be}-[(\Phi_{\be}\times \Phi_{\gm})\circ \Phi_{\al}](3)D^2(\xi_{1\al})\xi_{2\be}\\& &+(-1)^{\bar{\xi}_1+1}[(\Phi_{\be}\times \Phi_{\gm})\circ \Phi_{\al}](3)D(\xi_{1\al})D(\xi_{2\be})]+(-1)^{\bar{\xi}_3+1}[\la \Phi_{\be}\circ \Phi_{\gm}, \Phi_{\al}\ra D^7(\xi_{1\al})\xi_{2\be}\\&&+[(\Phi_{\be}\circ \Phi_{\gm})\cdot \Phi_{\al}](3)D^2(\xi_{1\al})\xi_{2\be}+[(\Phi_{\be}\circ \Phi_{\gm})\cdot \Phi_{\al}]D^4(\xi_{1\al})\xi_{2\be}\\& &+[(\Phi_{\be}\circ \Phi_{\gm})\times \Phi_{\al}](4)D(\xi_{1\al})\xi_{2\be}+[(\Phi_{\be}\circ \Phi_{\gm})\times \Phi_{\al}](2)D^3(\xi_{1\al})\xi_{2\be}\hspace{2cm}\end{eqnarray*}
\begin{eqnarray*}& &+[(\Phi_{\be}\circ\Phi_{\gm})\circ \Phi_{\al}](5)\xi_{1\al}\xi_{2\be}+[(\Phi_{\be}\circ \Phi_{\gm})\circ \Phi_{\al}](3)D^2(\xi_{1\al})\xi_{2\be}]\}\xi_{3\gm}.\hspace{2.4cm}(3.10)\end{eqnarray*}
Here we have viewed each term in (3.8-10) as an element in $\tilde{A}$ (cf. (2.23)).

For convenience, we call $D^{m_1}(\xi_{1\al})D^{m_2}(\xi_{2\be})\xi_{3\gm}$ a {\it monomial of index} $(0,m_1,m_2)$ and call $\Phi(n_1)D^{n_2}(\xi_{1\al})D^{n_3}(\xi_{2\be})\xi_{3\gm}$ a {\it monomial of index} $(n_1,n_2,n_3)$. We suppose that $H$ is Hamiltonian operator. Thus (2.33) holds. We substitute (3.8-10) into (2.33). By comparing the coefficients of the monomial of index (0,7,0) in (2.33), we have:
$$\la \Phi_{\al}\circ \Phi_{\be}, \Phi_{\gm}\ra+\la \Phi_{\be}\circ \Phi_{\gm}, \Phi_{\al}\ra=\la \Phi_{\be}\cdot \Phi_{\gm}, \Phi_{\al}\ra+\la \Phi_{\be}\times \Phi_{\gm}, \Phi_{\al}\ra.\eqno(3.11)$$
Moreover, by (3.5), (3.11) is equivalent to:
$$\la \Phi_{\al}\circ \Phi_{\be}, \Phi_{\gm}\ra+\la \Phi_{\gm}\circ \Phi_{\be}, \Phi_{\al}\ra=\la \Phi_{\al},\Phi_{\gm}\circ \Phi_{\be}\ra.\eqno(3.12)$$
Comparing the coefficients of the monomial of index (0,6,1) in (2.33), we obtain
$$\la \Phi_{\al}\times \Phi_{\be}, \Phi_{\gm}\ra+\la \Phi_{\be}\times \Phi_{\gm}, \Phi_{\al}\ra=\la \Phi_{\al}\circ \Phi_{\be}, \Phi_{\gm}\ra.\eqno(3.13)$$
The coefficients of the monomial of index (0,5,2) in (2.33) show:
$$\la \Phi_{\al}\cdot \Phi_{\be}, \Phi_{\gm}\ra+\la \Phi_{\be}\cdot \Phi_{\gm}, \Phi_{\al}\ra=3\la \Phi_{\al}\circ \Phi_{\be}, \Phi_{\gm}\ra.\eqno(3.14)$$
Examing the coefficients of the monomial of index (0,4,3) in (2.33), we have:
$$\la \Phi_{\al}\cdot \Phi_{\be}, \Phi_{\gm}\ra+3\la \Phi_{\al}\times \Phi_{\be}, \Phi_{\gm}\ra=3\la \Phi_{\al}\circ \Phi_{\be}, \Phi_{\gm}\ra.\eqno(3.15)$$
The coefficients of the monomial of index (0,3,4) in (2.33) imply
$$2\la \Phi_{\al}\cdot \Phi_{\be}, \Phi_{\gm}\ra=3\la \Phi_{\al}\circ \Phi_{\be}, \Phi_{\gm}\ra.\eqno(3.16)$$
Considering the coefficients of the monomial of index (0,2,5) in (2.33), we find:
$$\la \Phi_{\gm}\cdot \Phi_{\al}, \Phi_{\be}\ra+3\la \Phi_{\al}\circ \Phi_{\be}, \Phi_{\gm}\ra=2\la \Phi_{\al}\cdot \Phi_{\be}, \Phi_{\gm}\ra+3\la \Phi_{\al}\times \Phi_{\be}, \Phi_{\gm}\ra.\eqno(3.17)$$
Looking up the coefficients of the monomial of index (0,1,6) in (2.33), we have:
$$\la \Phi_{\gm}\times\Phi_{\al}, \Phi_{\be}\ra+\la \Phi_{\al}\circ \Phi_{\be}, \Phi_{\gm}\ra=\la \Phi_{\al}\cdot \Phi_{\be}, \Phi_{\gm}\ra.\eqno(3.18)$$
The following equation follows from  the coefficients of the monomial of index (0,0,7) in (2.33):
$$\la \Phi_{\gm}\circ\Phi_{\al}, \Phi_{\be}\ra+\la \Phi_{\al}\circ \Phi_{\be}, \Phi_{\gm}\ra=\la \Phi_{\al}\cdot \Phi_{\be}, \Phi_{\gm}\ra+\la \Phi_{\al}\times \Phi_{\be}, \Phi_{\gm}\ra.\eqno(3.19)$$

The coefficients of the monomial of index (5,0,0) in (2.33) tell us that
$$(\Phi_{\gm}\circ\Phi_{\al})\circ\Phi_{\be}+(\Phi_{\be}\circ \Phi_{\gm})\circ \Phi_{\al}=(\Phi_{\be}\cdot \Phi_{\gm})\circ \Phi_{\al}+(\Phi_{\be}\times \Phi_{\gm})\circ \Phi_{\al}.\eqno(3.20)$$
The coefficients of the monomial of index (4,1,0) in (2.33) give us the following equation:
$$(\Phi_{\gm}\times\Phi_{\al})\circ\Phi_{\be}+(\Phi_{\be}\circ \Phi_{\gm})\times \Phi_{\al}=(\Phi_{\be}\cdot \Phi_{\gm})\times \Phi_{\al}+(\Phi_{\be}\times \Phi_{\gm})\times \Phi_{\al}.\eqno(3.21)$$
Let us look at the coefficients of the monomial of index (4,0,1) in (2.33). We obtain:
$$(\Phi_{\gm}\circ\Phi_{\al})\times\Phi_{\be}=(\Phi_{\be}\times \Phi_{\gm})\circ \Phi_{\al}.\eqno(3.22)$$
Comparing the coefficients of the monomial of index (3,2,0) in (2.33), we get:
\begin{eqnarray*}& &(\Phi_{\gm}\cdot\Phi_{\al})\circ\Phi_{\be}+(\Phi_{\al}\circ \Phi_{\be})\cdot \Phi_{\gm}+(\Phi_{\al}\circ \Phi_{\be})\times \Phi_{\gm}-(\Phi_{\al}\circ \Phi_{\be})\circ \Phi_{\gm}\\& &+(\Phi_{\be}\circ \Phi_{\gm})\cdot \Phi_{\al}-(\Phi_{\be}\cdot \Phi_{\gm})\cdot \Phi_{\al}-(\Phi_{\be}\times \Phi_{\gm})\cdot \Phi_{\al} \\&=&(\Phi_{\be}\cdot \Phi_{\gm})\circ \Phi_{\al}+(\Phi_{\be}\times \Phi_{\gm})\circ \Phi_{\al}-(\Phi_{\be}\circ \Phi_{\gm})\circ \Phi_{\al}.\hspace{5cm}(3.23)\end{eqnarray*}
The coefficients of the monomial of index (3,1,1) in (2.33) tell us that
\begin{eqnarray*}& &(\Phi_{\gm}\times\Phi_{\al})\times\Phi_{\be}-(\Phi_{\gm}\times\Phi_{\al})\circ\Phi_{\be}+(\Phi_{\al}\times \Phi_{\be})\cdot \Phi_{\gm}+(\Phi_{\al}\times \Phi_{\be})\times \Phi_{\gm}\\&=&(\Phi_{\al}\times \Phi_{\be})\circ \Phi_{\gm}-(\Phi_{\be}\times \Phi_{\gm})\times \Phi_{\al} +(\Phi_{\be}\times \Phi_{\gm})\circ \Phi_{\al}.\hspace{4.5cm}(3.24)\end{eqnarray*}
Consulting the coefficients of the monomial of index (3,0,2) in (2.33), we find:
\begin{eqnarray*}& &(\Phi_{\gm}\circ\Phi_{\al})\cdot\Phi_{\be}+(\Phi_{\gm}\circ\Phi_{\al})\circ\Phi_{\be}-(\Phi_{\al}\cdot \Phi_{\be})\cdot \Phi_{\gm}-(\Phi_{\al}\cdot \Phi_{\be})\times \Phi_{\gm}\\&=& -(\Phi_{\al}\cdot \Phi_{\be})\circ \Phi_{\gm}+(\Phi_{\al}\times \Phi_{\be})\cdot \Phi_{\gm}+(\Phi_{\al}\times \Phi_{\be})\times \Phi_{\gm}-(\Phi_{\al}\times \Phi_{\be})\circ \Phi_{\gm}\\& &-(\Phi_{\al}\circ\Phi_{\be})\times \Phi_{\gm}-(\Phi_{\al}\circ \Phi_{\be})\cdot \Phi_{\gm}+(\Phi_{\al}\circ \Phi_{\be})\circ \Phi_{\gm}+(\Phi_{\be}\cdot\Phi_{\gm})\circ \Phi_{\al}.\hspace{1.6cm}(3.25)\end{eqnarray*}

The coefficients of the monomial of index (2,3,0) in (2.33) imply:
$$(\Phi_{\al}\circ\Phi_{\be})\times\Phi_{\gm}+(\Phi_{\be}\circ \Phi_{\gm})\times \Phi_{\al}=(\Phi_{\be}\cdot \Phi_{\gm})\times \Phi_{\al}+(\Phi_{\be}\times \Phi_{\gm})\times \Phi_{\al}.\eqno(3.26)$$
Extracting the coefficients of the monomial of index (2,2,1) in (2.33), we have:
$$(\Phi_{\gm}\cdot\Phi_{\al})\times\Phi_{\be}+(\Phi_{\al}\circ \Phi_{\be})\times \Phi_{\gm}-(\Phi_{\al}\times \Phi_{\be})\times \Phi_{\gm}=(\Phi_{\be}\times \Phi_{\gm})\times \Phi_{\al}+(\Phi_{\be}\times \Phi_{\gm})\cdot \Phi_{\al}.\eqno(3.27)$$
 The coefficients of the monomial of index (2,1,2) in (2.33) imply:
$$(\Phi_{\gm}\times\Phi_{\al})\cdot\Phi_{\be}+(\Phi_{\gm}\times\Phi_{\al})\times\Phi_{\be}+(\Phi_{\al}\circ \Phi_{\be})\times \Phi_{\gm}=(\Phi_{\al}\cdot \Phi_{\be})\times \Phi_{\gm}+(\Phi_{\be}\cdot \Phi_{\gm})\times \Phi_{\al}.\eqno(3.28)$$
By comparing the coefficients of the monomial of index (2,0,3) in (2.33), we have:
$$(\Phi_{\gm}\circ\Phi_{\al})\times\Phi_{\be}+(\Phi_{\al}\circ \Phi_{\be})\times \Phi_{\gm}=(\Phi_{\al}\cdot \Phi_{\be})\times \Phi_{\gm}+(\Phi_{\al}\times \Phi_{\be})\times \Phi_{\gm}.\eqno(3.29)$$

Comparing the coefficients of the monomial of index (1,4,0) in (2.33), we get:
$$(\Phi_{\al}\circ\Phi_{\be})\cdot\Phi_{\gm}+(\Phi_{\be}\circ \Phi_{\gm})\cdot \Phi_{\al}=(\Phi_{\be}\cdot \Phi_{\gm})\cdot \Phi_{\al}+(\Phi_{\be}\times \Phi_{\gm})\cdot \Phi_{\al}.\eqno(3.30)$$
The coefficients of the monomial of index (1,3,1) in (2.33) tell us that
$$(\Phi_{\al}\times\Phi_{\be})\cdot\Phi_{\gm}=\Phi_{\al}\cdot (\Phi_{\be}\times \Phi_{\gm}).\eqno(3.31)$$
The following equation follows from the coefficients of the monomial of index (1,2,2) in (2.33):
$$(\Phi_{\gm}\cdot\Phi_{\al})\cdot\Phi_{\be}-(\Phi_{\be}\cdot \Phi_{\gm})\cdot \Phi_{\al}+2(\Phi_{\al}\circ\Phi_{\be})\cdot\Phi_{\gm}=(\Phi_{\al}\cdot \Phi_{\be})\cdot \Phi_{\gm}+(\Phi_{\al}\times \Phi_{\be})\cdot \Phi_{\gm}.\eqno(3.32)$$
In terms of  the coefficients of the monomial of index (1,1,3) in (2.33), 
$$(\Phi_{\gm}\times\Phi_{\al})\cdot\Phi_{\be}=\Phi_{\gm}\cdot (\Phi_{\al}\times \Phi_{\be}).\eqno(3.33)$$
Consulting the coefficients of the monomial of index (1,0,4) in (2.33), we find:
$$(\Phi_{\gm}\circ\Phi_{\al})\cdot\Phi_{\be}+(\Phi_{\al}\circ \Phi_{\be})\cdot \Phi_{\gm}=(\Phi_{\al}\cdot \Phi_{\be})\cdot \Phi_{\gm}+(\Phi_{\al}\times \Phi_{\be})\cdot \Phi_{\gm}.\eqno(3.34)$$
Here we have always assumed that $\al,\be,\gm$ are three arbitrary elements of the index set $I$.

Next we shall do technical reductions. By (3.15) and (3.16), we have:
$${3\over 2}\la \Phi_{\al}\circ \Phi_{\be},\Phi_{\gm}\ra + 3\la \Phi_{\al}\times \Phi_{\be},\Phi_{\gm}\ra=3\la \Phi_{\al}\circ \Phi_{\be},\Phi_{\gm}\ra\Longrightarrow \la \Phi_{\al}\circ \Phi_{\be},\Phi_{\gm}\ra=2 \la \Phi_{\al}\times \Phi_{\be},\Phi_{\gm}\ra.\eqno(3.35)$$
Moreover, by (3.5) and (3.12), we can prove that (3.11-19) are equivalent to:
$$\la \Phi_{\al}\times \Phi_{\be},\Phi_{\gm}\ra=\la \Phi_{\al}, \Phi_{\be}\times\Phi_{\gm}\ra={1\over 2}\la \Phi_{\al}\circ \Phi_{\be},\Phi_{\gm}\ra={1\over 3}\la \Phi_{\al}\cdot \Phi_{\be},\Phi_{\gm}\ra.\eqno(3.36)$$

By (3.5), (3.20) is equivalent to:
$$(\Phi_{\gm}\circ \Phi_{\al})\circ \Phi_{\be}=(\Phi_{\gm}\circ \Phi_{\be})\circ \Phi_{\al}\eqno(3.37)$$
and (3.21) is equivalent to:
$$(\Phi_{\gm}\times \Phi_{\al})\circ \Phi_{\be}=(\Phi_{\gm}\circ \Phi_{\be})\times \Phi_{\al}.\eqno(3.38)$$
Note that (3.22) and (3.38) are equivalent. 

Again by (3.5), (3.23) is equivalent to:
$$(\Phi_{\gm}\cdot\Phi_{\al})\circ\Phi_{\be}+\Phi_{\gm}\circ (\Phi_{\al}\circ \Phi_{\be})=(\Phi_{\gm}\circ \Phi_{\be})\cdot \Phi_{\al}+(\Phi_{\gm}\circ \Phi_{\be})\circ \Phi_{\al},\eqno(3.39)$$
(3.24) is equivalent to:
$$(\Phi_{\gm}\times\Phi_{\al})\times\Phi_{\be}+(\Phi_{\be}\times \Phi_{\gm})\times \Phi_{\al}=(\Phi_{\al}\times \Phi_{\be})\circ \Phi_{\gm} +(\Phi_{\be}\times \Phi_{\gm})\circ \Phi_{\al}-\Phi_{\be}\circ(\Phi_{\gm}\times\Phi_{\al}),\eqno(3.40)$$
(3.25) is equivalent to:
$$(\Phi_{\gm}\circ\Phi_{\al})\cdot\Phi_{\be}+(\Phi_{\gm}\circ\Phi_{\al})\circ\Phi_{\be}=\Phi_{\gm}\circ(\Phi_{\be}\circ \Phi_{\al})+(\Phi_{\be}\cdot \Phi_{\gm})\circ \Phi_{\al},\eqno(3.41)$$
and (3.26), (2.29) are equivalent to (3.22). If we change the indices in (3.34) according to the cycle $\al\rightarrow\be\rightarrow \gm\rightarrow\al$, then we get (3.30). Similarly, (3.31) and (3.33) are equivalent. 
Furthermore,  (3.5), (3.32) and (3.34) imply:
$$(\Phi_{\gm}\cdot\Phi_{\al})\cdot\Phi_{\be}+(\Phi_{\al}\circ\Phi_{\be})\cdot\Phi_{\gm}=(\Phi_{\be}\cdot \Phi_{\gm})\cdot \Phi_{\al}+(\Phi_{\be}\circ \Phi_{\al})\cdot \Phi_{\gm},\eqno(3.42)$$
$$(\Phi_{\gm}\circ \Phi_{\al})\cdot \Phi_{\be}=(\Phi_{\be}\circ\Phi_{\al})\cdot \Phi_{\gm}.\eqno(3.43)$$

Our strategy to do further reduction is to get rid of ``$\cdot$'' in (3.27-28),
(3.31), (3.39) and (3.42-43) by (3.5) and (3.37-38). Note that (3.27) is equivalent to:
\begin{eqnarray*}& & (\Phi_{\gm}\circ\Phi_{\al})\times\Phi_{\be}+(\Phi_{\al}\circ\Phi_{\gm})\times\Phi_{\be}-(\Phi_{\gm}\times\Phi_{\al})\times\Phi_{\be}\\& &+(\Phi_{\al}\circ \Phi_{\be})\times \Phi_{\gm}-(\Phi_{\al}\times \Phi_{\be})\times \Phi_{\gm}\\&=&(\Phi_{\be}\times \Phi_{\gm})\times \Phi_{\al}+\Phi_{\al}\circ (\Phi_{\be}\times \Phi_{\gm})+ (\Phi_{\be}\times \Phi_{\gm})\circ \Phi_{\al}-(\Phi_{\be}\times \Phi_{\gm})\times \Phi_{\al},\hspace{1.1cm}(3.44)\end{eqnarray*}
which is equivalent to (3.40) by (3.38). 
Again using (3.5), (3.28) is equivalent to:
\begin{eqnarray*} & &(\Phi_{\gm}\times\Phi_{\al})\circ\Phi_{\be}+\Phi_{\be}\circ(\Phi_{\gm}\times\Phi_{\al})+(\Phi_{\al}\circ \Phi_{\be})\times \Phi_{\gm}-(\Phi_{\be}\circ \Phi_{\al})\times \Phi_{\gm}\\&=&-(\Phi_{\al}\times \Phi_{\be})\times \Phi_{\gm}
(\Phi_{\be}\circ \Phi_{\gm})\times \Phi_{\al}+(\Phi_{\gm}\circ \Phi_{\be})\times \Phi_{\al}-(\Phi_{\be}\times \Phi_{\gm})\times \Phi_{\al},\hspace{1.2cm}(3.45)\end{eqnarray*}
which is equivalent to (3.40) by (3.38). Furthermore, (3.31) is equivalent to:
\begin{eqnarray*} & & (\Phi_{\al}\times\Phi_{\be})\circ\Phi_{\gm}+\Phi_{\gm}\circ(\Phi_{\al}\times\Phi_{\be})-\Phi_{\al}\circ (\Phi_{\be}\times \Phi_{\gm})-(\Phi_{\be}\times \Phi_{\gm})\circ \Phi_{\al}\\&=&(\Phi_{\al}\times\Phi_{\be})\times\Phi_{\gm}-\Phi_{\al}\times(\Phi_{\be}\times \Phi_{\gm}).\hspace{7.6cm}(3.46)\end{eqnarray*}

Now (3.5) implies that (3.39) is equivalent to:
\begin{eqnarray*}& & (\Phi_{\gm}\circ\Phi_{\al})\circ\Phi_{\be}+(\Phi_{\al}\circ\Phi_{\gm})\circ\Phi_{\be}-(\Phi_{\gm}\times\Phi_{\al})\circ\Phi_{\be}+\Phi_{\gm}\circ (\Phi_{\al}\circ \Phi_{\be})\\&=&(\Phi_{\gm}\circ \Phi_{\be})\circ \Phi_{\al}+ \Phi_{\al}\circ (\Phi_{\gm}\circ \Phi_{\be})-(\Phi_{\gm}\circ \Phi_{\be})\times \Phi_{\al}+(\Phi_{\gm}\circ \Phi_{\be})\circ \Phi_{\al},\hspace{1.4cm}(3.47)\end{eqnarray*}
which by (3.37-38) is equivalent to:
$$(\Phi_{\al}\circ\Phi_{\gm})\circ\Phi_{\be}-\Phi_{\al}\circ (\Phi_{\gm}\circ \Phi_{\be})=(\Phi_{\gm}\circ \Phi_{\al})\circ \Phi_{\be}-\Phi_{\gm}\circ (\Phi_{\al}\circ \Phi_{\be}).\eqno(3.48)$$
Equations (3.37) and (3.48) shows that $(V,\circ)$ forms a Novikov algebra. 

Next by (3.5), (3.42) is equivalent to:
\begin{eqnarray*}& &(\Phi_{\gm}\circ\Phi_{\al})\cdot\Phi_{\be}+(\Phi_{\al}\circ\Phi_{\gm})\cdot\Phi_{\be}-(\Phi_{\gm}\times\Phi_{\al})\cdot\Phi_{\be}+(\Phi_{\al}\circ\Phi_{\be})\cdot\Phi_{\gm}\\&=&(\Phi_{\be}\circ \Phi_{\gm})\cdot \Phi_{\al}+(\Phi_{\gm}\circ \Phi_{\be})\cdot \Phi_{\al}-(\Phi_{\be}\times \Phi_{\gm})\cdot \Phi_{\al}+(\Phi_{\be}\circ \Phi_{\al})\cdot \Phi_{\gm},\hspace{2cm}(3.49)\end{eqnarray*}
which holds if (3.43) and (3.31) hold. Furthermore, (3.5) implies that (3.43) is equivalent to:
\begin{eqnarray*}&& (\Phi_{\gm}\circ \Phi_{\al})\circ \Phi_{\be}+\Phi_{\be}\circ(\Phi_{\gm}\circ \Phi_{\al})-(\Phi_{\gm}\circ \Phi_{\al})\times \Phi_{\be}\\&=&(\Phi_{\be}\circ\Phi_{\al})\circ \Phi_{\gm}+ \Phi_{\gm}\circ(\Phi_{\be}\circ\Phi_{\al})-(\Phi_{\be}\circ\Phi_{\al})\times \Phi_{\gm},\hspace{4.8cm}(3.50)\end{eqnarray*}
which holds if (3.37-38) and (3.48) hold. 

We summarize what we have proved as:
\psp

{\bf Theorem 3.1}. {\it A differential operator} $H$ {\it of the form (3.1) is a Hamiltonian operator if and only if} $(V,\circ,\times)$ {\it is an NX-bialgebra,}
$$u\cdot v=u\circ v+v\circ u-u\times v\qquad\;\;\mbox{for}\;\;u,v\in V,\eqno(3.51)$$
{\it  and} $\la\cdot,\cdot\ra$ {\it is a symmetric bilinear form satisfying:}
$$\la u\circ v,w\ra=\la u,v\circ w\ra =2\la u\times v,w\ra\qquad\;\;\mbox{for}\;\;u,v,w\in V.\eqno(3.52)$$
\psp

{\bf Example}. Let $({\cal A},\cdot,\circ )$ be a Novikov-Poisson algebra such that $({\cal A},\cdot)$ contains an identity element $1$ and
$$1\circ 1=2.\eqno(3.53)$$
We shall show now that $({\cal A},\cdot,\circ)$ is a NX-bialgebra. In fact,  we have
$$x\circ y=x\cdot \partial(y),\qquad\mbox{where}\;\;\partial(y)=1\circ y\eqno(3.54)$$
for $x,y\in {\cal A}$. Note by (1.16),
$$\partial(x\circ y)=\partial(x)\cdot y+x\cdot \partial(y)-2x\cdot y\qquad\mbox{for}\;\;x,y\in {\cal A}\eqno(3.55)$$
(cf. [X5]). Thus we have:
$$ (y\circ x)\cdot z+x\cdot (y\circ z)-y\circ (x\cdot z)=2y\cdot(x\cdot z)=(x\cdot y)\cdot z+x\cdot (y\cdot z)\eqno(3.56)$$
for $x,y,z\in {\cal A}$ by the commutativity and associativity of $({\cal A},\cdot).$ Furthermore,
\begin{eqnarray*}& & (x\cdot y)\circ z+z\circ (x\cdot y)-x\circ (y\cdot z)-(y\cdot z)\circ x\\&=& x\cdot y\cdot \partial(z)+z\cdot \partial(x)\cdot y+z\cdot x\cdot \partial(y)-x\cdot z\cdot \partial(y)\cdot z-x\cdot y\cdot \partial(z)-y\cdot z\cdot \partial(x)\\&=&0\\&=& (x\cdot y)\cdot z-x\cdot (y\cdot z)\hspace{10.1cm}(3.57)\end{eqnarray*}
for $x,y,z\in {\cal A}$. Hence $({\cal A},\cdot,\circ)$ is an NX-bialgebra. 

Next we shall give a concrete example. Let $({\cal A},\cdot)$ be the quotient algebra $\Bbb{R}[t]/(t^n)$ of the algebra $\Bbb{R}[t]$ of polynomials for a positive integer $n$. Denote by $e_j$ the image of $t^j$ in ${\cal A}$. We define the operation $\circ$ by
$$e_i\circ e_j=(j+2)e_{i+j}\qquad \mbox{for}\;\;0\leq i,j<n.\eqno(3.58)$$
Here we have used the convention that $e_l=0$ if $l\geq n$.
Then $({\cal A},\cdot,\circ)$ is a Novikov algebra satisfying (3.53) (cf. [X5]). Moreover, we define a bilinear form $\la\cdot,\cdot\ra$ on ${\cal A}$ by
$$\la e_i,e_j\ra=\delta_{i,0}\delta_{j,0}\qquad\mbox{for}\;\;0\leq i,j<n.\eqno(3.59)$$
Then $\la\cdot,\cdot\ra$ satisfies (3.52).

\section{Hamiltonian Superoperators and Fermionic Novikov Algebras}

Consider the following Hamiltonian operator $H$ of type 0:
$$-H^1_{\al,\be}=H^0_{\al,\be}=\sum_{\gm\in I}(a_{\al,\be}^{\gm}\Phi_{\gm}(2)+b_{\al,\be}^{\gm}\Phi_{\gm}D),\eqno(4.1)$$
where $a_{\al,\be}^{\gm},b_{\al,\be}^{\gm}\in \Bbb{R}$. Again we let $V$ be as in (3.6) and define operations: $\times,\circ:\;V\times V\rightarrow V$ by
$$\Phi_{\al}\circ \Phi_{\be}=\sum_{\gm\in I}a_{\al,\be}^{\gm}\Phi_r,\qquad \Phi_{\al}\times \Phi_{\be}=\sum_{\gm\in I}b_{\al,\be}^{\gm}\Phi_r\qquad\;\mbox{for}\;\;\al,\be\in I.\eqno(4.2)$$
\psp

{\bf Theorem 4.1}. {\it A differential operator of the form (4.1) is a Hamiltonian operator if and only if} $(V,\circ)$ {\it is a fermionic Novikov algebra and}
$$u\times v=v\circ u-u\circ v\qquad\qquad\mbox{\it for}\;\;u,v\in V.\eqno(4.3)$$
\psp

{\it Proof.} By (2.31), the super skew-symmetry of $H$ is equivalent to:
\begin{eqnarray*}\sum_{\gm\in I}(a_{\al,\be}^{\gm}\Phi_{\gm}(2)+b_{\al,\be}^{\gm}D\circ \Phi_{\gm})&=&\sum_{\gm\in I}((a_{\al,\be}^{\gm}+b_{\al,\be}^{\gm})\Phi_{\gm}(2)-b_{\al,\be}^{\gm}\Phi_{\gm}D)\\&=&\sum_{\gm\in I}(a_{\be,\al}^{\gm}\Phi_{\gm}(2)+b_{\be,\al}^{\gm}\Phi_{\gm}D)\hspace{4.6cm}(4.4)\end{eqnarray*}
for $\al,\be\in I$, which is equivalent to (4.3). Next we shall find the exact formula for each term in (3.33).  For any $\bar{\xi}_1,\bar{\xi}_2,\bar{\xi}_3\in \Omega$, we have:
\begin{eqnarray*}& &\bar{\xi}_3((D_H\bar{\xi}_1)H\bar{\xi}_2)\\
&=& (-1)^{\bar{\xi}_1+\bar{\xi}_2}\sum_{\al,\be,\gm,\lmd,\mu\in I}[(-1)^{\bar{\xi}_1}a_{\gm,\al}^{\lmd}\xi_{1\al}D+b_{\gm,\al}^{\lmd}D(\xi_{1\al})][a_{\lmd,\be}^{\mu}\Phi_{\mu}(2)\xi_{2\be}+b_{\lmd,\be}^{\mu}\Phi_{\mu}D(\xi_{2\be})]\xi_{3\gm}\\&=& \sum_{\al,\be,\gm\in I}\{(-1)^{\bar{\xi}_1+\bar{\xi}_2+1}[(\Phi_{\gm}\circ \Phi_{\al})\circ \Phi_{\be}](3)\xi_{1\al}\xi_{2\be}+(-1)^{\bar{\xi}_2}[(\Phi_{\gm}\circ \Phi_{\al})\circ \Phi_{\be}](2)\xi_{1\al}D(\xi_{2\be})\\&&+(-1)^{\bar{\xi}_2}[(\Phi_{\gm}\circ \Phi_{\al})\times \Phi_{\be}](2)\xi_{1\al}D(\xi_{2\be})+(-1)^{\bar{\xi}_1+\bar{\xi}_2}[(\Phi_{\gm}\circ \Phi_{\al})\times \Phi_{\be}]\xi_{1\al}D^2(\xi_{2\be})\\& &+(-1)^{\bar{\xi}_1+\bar{\xi}_2}[(\Phi_{\gm}\times \Phi_{\al})\circ \Phi_{\be}](2)D(\xi_{1\al})\xi_{2\be}\\& &+(-1)^{\bar{\xi}_2}[(\Phi_{\gm}\times \Phi_{\al})\times \Phi_{\be}]D(\xi_{1\al})D(\xi_{2\be})\}\xi_{3\gm},\hspace{6.4cm}(4.5)\end{eqnarray*}
\begin{eqnarray*}& &(-1)^{(\bar{\xi}_1+1)(\bar{\xi}_2+\bar{\xi}_3)}\bar{\xi}_1((D_H\bar{\xi}_2)H\bar{\xi}_3)\\&=&
\sum_{\al,\be,\gm\in I}\{(-1)^{\bar{\xi}_2+\bar{\xi}_3+1}[(\Phi_{\al}\circ \Phi_{\be})\circ \Phi_{\gm}](3)\xi_{1\al}\xi_{2\be}+(-1)^{\bar{\xi}_2+\bar{\xi}_3}[[(\Phi_{\al}\circ \Phi_{\be})\circ \Phi_{\gm}](3)\xi_{1\al}\xi_{2\be}\\& &+[(\Phi_{\al}\circ \Phi_{\be})\circ \Phi_{\gm}](2)D(\xi_{1\al})\xi_{2\be}+(-1)^{\bar{\xi}_1+1}[(\Phi_{\al}\circ \Phi_{\be})\circ \Phi_{\gm}](2)\xi_{1\al}D(\xi_{2\be})]
\\&&+(-1)^{\bar{\xi}_2+\bar{\xi}_3}[[(\Phi_{\al}\circ \Phi_{\be})\times \Phi_{\gm}](3)\xi_{1\al}\xi_{2\be}+[(\Phi_{\al}\circ \Phi_{\be})\times \Phi_{\gm}](2)D(\xi_{1\al})\xi_{2\be}\\& &+(-1)^{\bar{\xi}_1+1}[(\Phi_{\al}\circ \Phi_{\be})\times \Phi_{\gm}](2)\xi_{1\al}D(\xi_{2\be})]+(-1)^{\bar{\xi}_2+\bar{\xi}_3+1}[[(\Phi_{\al}\circ \Phi_{\be})\times \Phi_{\gm}](3)\xi_{1\al}\xi_{2\be}\\& &+[(\Phi_{\al}\circ \Phi_{\be})\times \Phi_{\gm}]D^2(\xi_{1\al})\xi_{2\be}+[(\Phi_{\al}\circ \Phi_{\be})\times \Phi_{\gm}]\xi_{1\al}D^2(\xi_{2\be})]
\\& &+(-1)^{\bar{\xi}_1+\bar{\xi}_3+1}[(\Phi_{\al}\times \Phi_{\be})\circ \Phi_{\gm}](2)\xi_{1\al}D(\xi_{2\be})+(-1)^{\bar{\xi}_1+\bar{\xi}_2+\bar{\xi}_3+1}[[(\Phi_{\al}\times \Phi_{\be})\circ \Phi_{\gm}](2)\xi_{1\al}D(\xi_{2\be})\\& &-[(\Phi_{\al}\times \Phi_{\be})\circ \Phi_{\gm}]D(\xi_{1\al})D(\xi_{2\be})+(-1)^{\bar{\xi}_1}[(\Phi_{\al}\times \Phi_{\be})\times \Phi_{\gm}]\xi_{1\al}D^2(\xi_{2\be})]\}\xi_{3\gm},\hspace{0.8cm}(4.6)\end{eqnarray*}
\begin{eqnarray*}& &(-1)^{(\bar{\xi}_3+1)(\bar{\xi}_1+\bar{\xi}_2)}\bar{\xi}_2((D_H\bar{\xi}_3)H\bar{\xi}_1)\\&=& \sum_{\al,\be,\gm\in I}\{(-1)^{\bar{\xi}_1+\bar{\xi}_3+1}[(\Phi_{\be}\circ \Phi_{\gm})\circ \Phi_{\al}](3)\xi_{1\al}\xi_{2\be}+(-1)^{\bar{\xi}_1+\bar{\xi}_3+1}[(\Phi_{\be}\circ \Phi_{\gm})\circ \Phi_{\al}](2)D(\xi_{1\al})\xi_{2\be}\\&&+(-1)^{\bar{\xi}_1+\bar{\xi}_3+1}[(\Phi_{\be}\circ \Phi_{\gm})\times \Phi_{\al}](2)D(\xi_{1\al})\xi_{2\be}+(-1)^{\bar{\xi}_3+\bar{\xi}_1}[(\Phi_{\be}\circ \Phi_{\gm})\times \Phi_{\al}]D^2(\xi_{1\al})\xi_{2\be}\\& &+(-1)^{\bar{\xi}_3+\bar{\xi}_1+1}[[(\Phi_{\be}\times \Phi_{\gm})\circ \Phi_{\al}](3)\xi_{1\al}\xi_{2\be}+[(\Phi_{\be}\times \Phi_{\gm})\circ \Phi_{\al}](2)D(\xi_{1\al})\xi_{2\be}\\&&(-1)^{\bar{\xi}_1+1}[(\Phi_{\be}\times \Phi_{\gm})\circ \Phi_{\al}](2)\xi_{1\al}D(\xi_{2\be})]+(-1)^{\bar{\xi}_1+\bar{\xi}_3+1}[[(\Phi_{\be}\times \Phi_{\gm})\times\Phi_{\al}](2)D(\xi_{1\al})\xi_{2\al}\\&& -[(\Phi_{\be}\times \Phi_{\gm})\times \Phi_{\al}]D^2(\xi_{1\al})\xi_{2\al}+(-1)^{\bar{\xi}_1+1}[(\Phi_{\be}\times \Phi_{\gm})\times \Phi_{\al}]D(\xi_{1\al})D(\xi_{2\al})]\}\xi_{3\gm}.\hspace{0.5cm}(4.7)\end{eqnarray*}
We assume that $H$ is a Hamiltonian operator. Thus (2.33) holds. We substitute (4.5-7) into (2.33). We define the monomial index as in Section 3. In the following, we always assume that $\al,\be,\gm$ are arbitrary elements of $I$. By comparing the coefficients of the monomial of index (3,0,0) in (2.33), we have:
$$(\Phi_{\gm}\circ \Phi_{\al})\circ \Phi_{\be}+(\Phi_{\be}\circ \Phi_{\gm})\circ \Phi_{\al}+(\Phi_{\be}\times \Phi_{\gm})\circ \Phi_{\al}=0,\eqno(4.8)$$
which by (4.3) is equivalent to:
$$(\Phi_{\gm}\circ \Phi_{\al})\circ \Phi_{\be}=-(\Phi_{\gm}\circ \Phi_{\be})\circ \Phi_{\al}.\eqno(4.9)$$
The coefficients of the monomial of index (2,1,0) in (2.33) imply:
\begin{eqnarray*}& &(\Phi_{\gm}\times \Phi_{\al})\circ \Phi_{\be}+(\Phi_{\al}\circ \Phi_{\be})\circ \Phi_{\gm}+(\Phi_{\al}\circ \Phi_{\be})\times \Phi_{\gm}-(\Phi_{\be}\circ \Phi_{\gm})\circ \Phi_{\al}\\&=&(\Phi_{\be}\circ \Phi_{\gm})\times \Phi_{\al}+(\Phi_{\be}\times \Phi_{\gm})\circ \Phi_{\al}+(\Phi_{\be}\times \Phi_{\gm})\times \Phi_{\al},\hspace{4.5cm}(4.10)\end{eqnarray*}
which by (4.3) is equivalent to:
$$(\Phi_{\gm}\times \Phi_{\al})\circ \Phi_{\be}+ \Phi_{\gm}\circ(\Phi_{\al}\circ \Phi_{\be})-\Phi_{\al}\circ (\Phi_{\be}\circ \Phi_{\gm})-\Phi_{\al}\circ (\Phi_{\be}\times \Phi_{\gm})=0.\eqno(4.11)$$
Again by (4.3), (4.11) is equivalent to:
$$(\Phi_{\al}\circ \Phi_{\gm})\circ \Phi_{\be}-\Phi_{\al}\circ (\Phi_{\gm}\circ \Phi_{\be})=(\Phi_{\gm}\circ \Phi_{\al})\circ \Phi_{\be}- \Phi_{\gm}\circ(\Phi_{\al}\circ \Phi_{\be}).\eqno(4.12)$$
Note that (4.9) and (4.12) imply that $(V,\circ)$ is a fermionic Novikov algebra. Consulting the coefficients of the monomial of index (2,0,1) in (2.33), we get:
\begin{eqnarray*}& &(\Phi_{\gm}\circ \Phi_{\al})\circ \Phi_{\be}+(\Phi_{\gm}\circ \Phi_{\al})\times \Phi_{\be}-(\Phi_{\al}\circ \Phi_{\be})\circ \Phi_{\gm}-(\Phi_{\al}\circ \Phi_{\be})\times \Phi_{\gm}\\&=&(\Phi_{\al}\times \Phi_{\gm})\circ \Phi_{\gm}+(\Phi_{\al}\times \Phi_{\be})\times \Phi_{\gm}-(\Phi_{\be}\times \Phi_{\gm})\circ\Phi_{\al},\hspace{4.5cm}(4.13)\end{eqnarray*}
which is (4.10) if we change the indices according to the cycle $\al\rightarrow \be\rightarrow\gm\rightarrow\al$. 

Examing the coefficients of the monomial of index (1,2,0) in (2.33), we have:
$$(\Phi_{\al}\circ \Phi_{\be})\times \Phi_{\gm}-(\Phi_{\be}\circ \Phi_{\gm})\times \Phi_{\al}-(\Phi_{\be}\times \Phi_{\gm})\times \Phi_{\al}=0,\eqno(4.14)$$
which by (4.3) is equivalent to:
$$(\Phi_{\al}\circ \Phi_{\be})\times \Phi_{\gm}-(\Phi_{\gm}\circ \Phi_{\be})\times \Phi_{\al}=0.\eqno(4.15)$$
Again by (4.3), (4.15) is equivalent to:
$$\Phi_{\gm}\circ(\Phi_{\al}\circ \Phi_{\be})-(\Phi_{\al}\circ \Phi_{\be})\circ \Phi_{\gm}- \Phi_{\al}\circ(\Phi_{\gm}\circ \Phi_{\be})+(\Phi_{\gm}\circ \Phi_{\be})\circ \Phi_{\al}=0,\eqno(4.16)$$
which is equivalent to (4.12) by (4.9). The coefficients of the monomial of index (1,0,2) in (2.33) tell us that
$$(\Phi_{\gm}\circ \Phi_{\al})\times \Phi_{\be}-(\Phi_{\al}\circ \Phi_{\be})\times \Phi_{\gm}-(\Phi_{\al}\times \Phi_{\be})\times \Phi_{\gm}=0,\eqno(4.17)$$
which is (4.14) if we change indices according to the cycle $\al\rightarrow \be\rightarrow\gm\rightarrow\al$.
Finally,  checking the coefficients of the monomial of index (1,1,1) in (2.33), we obtain:
$$(\Phi_{\gm}\times \Phi_{\al})\times \Phi_{\be}+(\Phi_{\al}\times \Phi_{\be})\times \Phi_{\gm}+(\Phi_{\be}\times \Phi_{\gm})\times \Phi_{\al}=0,\eqno(4.18)$$
which by (4.3) is equivalent to:
\begin{eqnarray*}& &(\Phi_{\al}\circ \Phi_{\gm})\times \Phi_{\be}-(\Phi_{\gm}\circ \Phi_{\al})\times \Phi_{\be}+(\Phi_{\be}\circ\Phi_{\al})\times \Phi_{\gm}\\&=&(\Phi_{\al}\circ \Phi_{\be})\times \Phi_{\gm}-(\Phi_{\gm}\circ \Phi_{\be})\times \Phi_{\al}+(\Phi_{\be}\circ \Phi_{\gm})\times \Phi_{\al},\hspace{4.5cm}(4.19)\end{eqnarray*}
which holds if (4.15) is satisfied. This shows that $(V,\times)$ is a Lie algebra. From the above arguments, one can see that we have proved that a matrix differential operator $H$ of the form (4.1) is a Hamiltonian operator if and only if (4.3), (4.9) and (4.12) are satisfied.\hspace{1cm}$\Box$
\psp

{\bf Example}. It is not that easy to construct nontrivial fermionic Novikov algebras. Let $V$ be a vector space with a basis $\{e_1,e_2,e_3,e_4\}$ and $(E(V),\cdot)$ be the exterior algebra generated by $V$. Then $E(V)$ is 16-dimensional. Set
$$v_1=e_2\cdot e_3\cdot e_4,\;\;v_2=e_1\cdot e_3\cdot e_4,\;\;v_3=e_1\cdot e_2\cdot e_4,\;\;v_4=e_1\cdot e_2\cdot e_3,\eqno(4.20)$$
$$v_0=\sum_{i<j}c_{i,j}e_i\cdot e_j,\qquad v_5=e_1\cdot e_2\cdot e_3\cdot e_4,\eqno(4.21)$$
where $c_{i,j}\in \Bbb{R}$ are constants. We define
$${\cal A}=\sum_{i=0}^5\Bbb{R}v_i\eqno(4.22)$$
and define the operation on ${\cal A}$ by:
$$v\circ v_0=v\circ v_5=0,\qquad v\circ v_i=v\cdot e_i\qquad\qquad\mbox{for}\;\;v\in {\cal A},\;i=1,2,3,4.\eqno(4.23)$$
Then the operation $\circ$ satisfies (4.9).  

Let us prove (4.12), that is
$$(v_i\circ v_j)\circ v_k-v_i\circ (v_j\circ v_k)=(v_j\circ v_i)\circ v_k-v_j\circ (v_i\circ v_k)\qquad\mbox{for}\;\;i,j,k=0,1,...,5.\eqno(4.24)$$
Notice that (4.24) holds obviously  if one of the following conditions is satisfied : (a)$i=j$;  (b) $k=0,5$; (c) $i=5$; (d) $j=5$; (e) $i,j\in \{1,2,3,4\}$. Moreover, by symmetry, we only need to prove it when $i=0,\;j=1$ and $k=1$ or 2. 
\begin{eqnarray*}\hspace{2cm}& &(v_0\circ v_1)\circ v_1-v_0\circ (v_1\circ v_1)\\&=&v_0\cdot e_1\cdot e_1-v_0\circ (v_1\cdot e_1)\\&=& v_0\circ v_5\\&=&0,\hspace{11.6cm}(4.25)\end{eqnarray*}
\begin{eqnarray*}\hspace{2cm}& &(v_1\circ v_0)\circ v_1-v_1\circ (v_0\circ v_1)\\&=&-v_1\circ (v_0\cdot e_1)\\&=& -v_1\circ (c_{2,3}v_4+c_{2,4}v_3+c_{3,4}v_2)\\&=&-(e_2\cdot e_3\cdot e_4)\cdot (c_{2,3}e_4+c_{2,4}e_3+c_{3,4}e_2)\\&=&0;\hspace{11.6cm}(4.26)\end{eqnarray*}
\begin{eqnarray*}\hspace{2cm}& &(v_0\circ v_1)\circ v_2-v_0\circ (v_1\circ v_2)\\&=&v_0\cdot e_1\cdot e_2-v_0\circ (v_1\cdot e_2)\\&=&c_{3,4}v_5,\hspace{10.9cm}(4.27)\end{eqnarray*}
\begin{eqnarray*}\hspace{2cm}& &(v_1\circ v_0)\circ v_2-v_1\circ (v_0\circ v_2)\\&=&-v_1\circ (v_0\cdot e_2)\\&=& -v_1\circ (-c_{1,3}v_4-c_{1,4}v_3+c_{3,4}v_1)\\&=&-(e_2\cdot e_3\cdot e_4)\cdot (-c_{1,3}e_4-c_{1,4}e_3+c_{3,4}e_1)\\&=&c_{3,4}v_5.\hspace{10.9cm}(4.28)\end{eqnarray*}
Thus we prove that the algebra $({\cal A},\circ)$ defined in (4.21-23) is a fermionic Novikov algebra.

\section{Induced Lie Superalgebras}

In this section, we shall prove that a type-1 Hamiltonian operator $H$ of the form
$$H^1_{\al,\be}=H^0_{\al,\be}=\sum_{\gm\in I}[\sum_{m=0}^Na_{\al,\be,\gm}^m\Phi_{\gm}(2(N-m)+1)D^{2m}+\sum_{n=0}^{N-1}b_{\al,\be,\gm}^n\Phi_{\gm}(2(N-n))D^{2n+1}]\eqno(5.1)$$
induces a Lie superalgebra.

In the rest of this section, we denote by $\theta_i$ anticommuative formal variables and by $z_i$ commutative formal variables for $i=1,2,3$, that is,
$$\theta_i\theta_j=-\theta_j\theta_i,\qquad z_i\theta_j=\theta_j z_i,\qquad z_iz_j=z_jz_i\qquad\mbox{for}\;\;i,j=1,2,3.\eqno(5.2)$$
We let
$$\delta\left({z_i\over z_j}\right)=\sum_{m\in \Bbb{Z}}{z_i^m\over z_j^m},\qquad \Delta_{i,j}=(\theta_i-\theta_j)\delta\left({z_i\over z_j}\right).\eqno(5.3)$$
Note that
$$\Delta_{i,j}=-\Delta_{j,i}.\eqno(5.4)$$
Let
$$f(\theta,z)=f_0(z)+\theta f_1(z),\qquad\mbox{for}\;\;f_i(z)\in \Bbb{R}[z,z^{-1}].\eqno(5.5)$$
\psp

{\bf Lemma 5.1}. {\it We have}:
$$f(\theta_1,z_1)\Delta_{1,2}=f(\theta_2,z_2)\Delta_{1,2}.\eqno(5.6)$$

{\it Proof.}
\begin{eqnarray*}\hspace{1cm}&&f(\theta_1,z_1)\Delta_{1,2}\\&=&(f_0(z_1)+\theta_1f_1(z_1))(\theta_1-\theta_2)\delta\left({z_1\over z_2}\right)\\&=&(\theta_1-\theta_2)f_0(z_1)\delta\left({z_1\over z_2}\right)-\theta_1\theta_2f_1(z_1)\delta\left({z_1\over z_2}\right)\\&=& (\theta_1-\theta_2)f_0(z_2)\delta\left({z_1\over z_2}\right)+\theta_2\theta_1f_1(z_2)\delta\left({z_1\over z_2}\right)\\&=&(f_0(z_2)+\theta_2f_1(z_2))(\theta_1-\theta_2)\delta\left({z_1\over z_2}\right)\\&=&f(\theta_2,z_2)\Delta_{1,2}.\qquad\qquad\Box\hspace{8.8cm}(5.7)\end{eqnarray*}

Let $L$ be a vector space with a basis $\{\phi_{\al}(n)\mid \al\in I,n\in \Bbb{Z}/2\}$. We denote
$$\phi_{\al}(\theta,z)=\sum_{n\in\Bbb{Z}}\phi_{\al}(n)z^{-n-N-1}\theta+\sum_{n\in\Bbb{Z}}\phi_{\al}\left(n+{1\over 2}\right)z^{-n-N-1}=\phi_{\al}^0(z)\theta+\phi_{\al}^1(z)\eqno(5.8)$$
for $\al\in I.$ Our notions are motivated by the theory of vertex operator algebras (e.g., cf. [FLM]). In the rest of this section, we always assume that
$$\theta_i \phi_{\al}(n)=\phi_{\al}(n)\theta_i,\;\;\;\theta_i \phi_{\al}\left(n+{1\over 2}\right)=-\phi_{\al}\left(n+{1\over 2}\right)\theta_i;\eqno(5.9)$$
$$z_i \phi_{\al}(n)=\phi_{\al}(n)z_i,\;\;\;z_i \phi_{\al}\left(n+{1\over 2}\right)=\phi_{\al}\left(n+{1\over 2}\right)z_i\eqno(5.10)$$
for $\al\in I,n\in \Bbb{Z}$ and $i=1,2,3.$ 

We also use the notions:
$$D_i=\theta_i\partial_{z_i}+\partial_{\theta_i}\qquad\mbox{for}\;\;i=1,2,3;\eqno(5.11)$$
By induction on $n\in \Bbb{N}$, we can prove:
\psp

{\bf Lemma 5.2}. 
$$z_2^{-1}D_1^{2n}\Delta_{1,2}=(-1)^nz_1^{-1}D_2^{2n}\Delta_{1,2},\qquad z_2^{-1}D_1^{2n+1}\Delta_{1,2}=(-1)^{n+1}z_1^{-1}D_2^{2n+1}\Delta_{1,2}\eqno(5.12)$$
{\it for} $n\in \Bbb{N}.$
\psp

Now we define the operation $[\cdot,\cdot]$ on $L$ by
\begin{eqnarray*}[\phi_{\al}(\theta_1,z_1),\phi_{\be}(\theta_2,z_2)]&=&z_2^{-1}\sum_{\gm\in I}\{\sum_{m=0}^Na_{\al,\be,\gm}^mD_1^{2(N-m)}(\phi_{\gm}(\theta_1,z_1))D_1^{2m}(\Delta_{1,2})\\& &
+\sum_{n=0}^{N-1}b_{\al,\be,\gm}^nD_1^{2(N-n)-1}(\phi_{\gm}(\theta_1,z_1))D_1^{2n+1}(\Delta_{1,2})\}\hspace{2cm}(5.13)\end{eqnarray*}
for $\al,\be \in I$.
\psp

{\bf Theorem 5.3}. {\it The algebra} $(L,[\cdot,\cdot])$ {\it forms a Lie supealgebra with the grading:}
$$L_0=\sum_{\al\in I}\sum_{n\in \Bbb{Z}}\Bbb{R}\phi_{\al}(n),\qquad L_1=\sum_{\al\in I}\sum_{n\in \Bbb{Z}}\Bbb{R}\phi_{\al}\left(n+{1\over 2}\right).\eqno(5.14)$$

{\it Proof}. First, we have:
\begin{eqnarray*}\hspace{1cm}& &[\phi_{\al}(\theta_1,z_1),\phi_{\be}(\theta_2,z_2)]\\&\stackrel{\tiny (5.12)}{=}&z_1^{-1}\sum_{\gm\in I}\{\sum_{m=0}^N(-1)^ma_{\al,\be,\gm}^mD_1^{2(N-m)}(\phi_{\gm}(\theta_1,z_1))D_2^{2m}(\Delta_{1,2})\\& &
+\sum_{n=0}^{N-1}(-1)^{n+1}b_{\al,\be,\gm}^nD_1^{2(N-n)-1}(\phi_{\gm}(\theta_1,z_1))D_2^{2n+1}(\Delta_{1,2})\}\\&=&z_1^{-1}\sum_{\gm\in I}\{\sum_{m=0}^N(-1)^ma_{\al,\be,\gm}^mD_2^{2m}[D_1^{2(N-m)}(\phi_{\gm}(\theta_1,z_1))\Delta_{1,2}]\\& &
+\sum_{n=0}^{N-1}(-1)^{n+1}b_{\al,\be,\gm}^nD_2^{2n+1}[D_1^{2(N-n)-1}(\phi_{\gm}(\theta_1,z_1))\Delta_{1,2}]\}\\&\stackrel{\tiny (5.6)}{=}&z_1^{-1}\sum_{\gm\in I}\{\sum_{m=0}^N(-1)^ma_{\al,\be,\gm}^mD_2^{2m}[D_2^{2(N-m)}(\phi_{\gm}(\theta_2,z_2))\Delta_{1,2}]\hspace{5cm}\end{eqnarray*}
\begin{eqnarray*}& &
+\sum_{n=0}^{N-1}(-1)^{n+1}b_{\al,\be,\gm}^nD_2^{2n+1}[D_2^{2(N-n)-1}(\phi_{\gm}(\theta_2,z_2))\Delta_{1,2}]\}\\&=&z_1^{-1}\sum_{\gm\in I}\{\sum_{m=0}^N(-1)^{m+1}a_{\al,\be,\gm}^mD_2^{2m}[D_2^{2(N-m)}(\phi_{\gm}(\theta_2,z_2))\Delta_{2,1}]\\& &+\sum_{n=0}^{N-1}(-1)^nb_{\al,\be,\gm}^nD_2^{2n+1}[D_2^{2(N-n)-1}(\phi_{\gm}(\theta_2,z_2))\Delta_{2,1}]\}\hspace{4.5cm}(5.15)\end{eqnarray*}
for $\al,\be\in I$. Therefore, the super skew-symmetry of $H$ and (2.31) imply
\begin{eqnarray*}& &[\phi_{\al}^0(z_1),\phi_{\be}^0(z_2)]\theta_1\theta_2-[\phi_{\al}^0(z_1),\phi_{\be}^1(z_2)]\theta_1+[\phi_{\al}^1(z_1),\phi_{\be}^0(z_2)]\theta_2+[\phi_{\al}^1(z_1),\phi_{\be}^1(z_2)]\\&=&
[\phi_{\al}(\theta_1,z_1),\phi_{\be}(\theta_2,z_2)]\\&=&[\phi_{\be}(\theta_2,z_2),\phi_{\al}(\theta_1,z_1)]\\&=&[\phi_{\be}^1(z_2),\phi_{\al}^0(z_1)]\theta_2\theta_1-[\phi_{\be}^0(z_2),\phi_{\al}^1(z_1)]\theta_2\\& &+[\phi_{\be}^1(z_2),\phi_{\al}^0(z_1)]\theta_1+[\phi_{\be}^1(z_2),\phi_{\al}^1(z_1)]\hspace{7.4cm}(5.16)\end{eqnarray*}
for $\al,\be\in I$, which implies the skew-symmetry:
$$[\phi_{\al}^i(z_1),\phi^j_{\be}(z_2)]=-(-1)^{ij}[\phi^j_{\be}(z_2),\phi_{\al}^i(z_1)]\qquad\mbox{for}\;\;\al,\be\in I;\;i,j\in \Bbb{Z}_2.\eqno(5.17)$$

In the rest of this section, we assume that $\al,\be,\gm$ are arbitrary elements of $I$. 
Note that
\begin{eqnarray*}&& [[\phi_{\al}(\theta_1,z_1),\phi_{\be}(\theta_2,z_2)] ,\phi_{\gm}(\theta_3,z_3)]\\&=& z_2^{-1}\sum_{\lmd\in I}\{\sum_{m=0}^Na_{\al,\be,\lmd}^m[D_1^{2(N-m)}(\phi_{\lmd}(\theta_1,z_1))D_1^{2m}(\Delta_{1,2}),\phi_{\gm}(\theta_3,z_3)]\\& &
+\sum_{m=0}^{N-1}b_{\al,\be,\lmd}^m[D_1^{2(N-m)-1}(\phi_{\lmd}(\theta_1,z_1))D_1^{2m+1}(\Delta_{1,2}),\phi_{\lmd}(\theta_3,z_3)]\}\\&=& z_2^{-1}\sum_{\lmd\in I}\{\sum_{m=0}^N-a_{\al,\be,\lmd}^mD_1^{2(N-m)}[\phi_{\lmd}(\theta_1,z_1),\phi_{\gm}(\theta_3,z_3)]D_1^{2m}(\Delta_{1,2})\\& &
+\sum_{m=0}^{N-1}b_{\al,\be,\lmd}^mD_1^{2(N-m)-1}[\phi_{\lmd}(\theta_1,z_1),\phi_{\lmd}(\theta_3,z_3)]D_1^{2m+1}(\Delta_{1,2})\}\\&=& z_2^{-1}z_3^{-1}\sum_{\lmd,\mu\in I}\{\sum_{m=0}^N\sum_{n=0}^N-a_{\al,\be,\lmd}^ma_{\lmd,\gm,\mu}^nD_1^{2(N-m)}[D_1^{2(N-n)}(\phi_{\mu}(\theta_1,z_1))D_1^{2n}(\Delta_{1,3})]D_1^{2m}(\Delta_{1,2})\\& &+\sum_{m=0}^N\sum_{n=0}^{N-1}-a_{\al,\be,\lmd}^mb_{\lmd,\gm,\mu}^nD_1^{2(N-m)}[D_1^{2(N-n)-1}(\phi_{\mu}(\theta_1,z_1))D_1^{2n+1}(\Delta_{1,3})]D_1^{2m}(\Delta_{1,2})\\& &+\sum_{m=0}^{N-1}\sum_{n=0}^Nb_{\al,\be,\lmd}^ma_{\lmd,\gm,\mu}^nD_1^{2(N-m)-1}[D_1^{2(N-n)}(\phi_{\mu}(\theta_1,z_1))D_1^{2n}(\Delta_{1,3})]D_1^{2m+1}(\Delta_{1,2})
\hspace{2cm}\end{eqnarray*}
\begin{eqnarray*}&&+\sum_{m=0}^{N-1}\sum_{n=0}^{N-1}b_{\al,\be,\lmd}^mb_{\lmd,\gm,\mu}^nD_1^{2(N-m)-1}[D_1^{2(N-n)-1}(\phi_{\mu}(\theta_1,z_1))D_1^{2n+1}(\Delta_{1,3})]D_1^{2m+1}(\Delta_{1,2})\},\hspace{0.6cm}(5.18)\end{eqnarray*}

\begin{eqnarray*}&& [[\phi_{\be}(\theta_2,z_2),\phi_{\gm}(\theta_3,z_3)] ,\phi_{\al}(\theta_1,z_1)]\\&=&z_2^{-1}z_3^{-1}\sum_{\lmd,\mu\in I}\{\sum_{m=0}^N\sum_{n=0}^N(-1)^{N+m+n+1}a_{\be,\gm,\lmd}^ma_{\lmd,\al,\mu}^nD_1^{2n}[D_1^{2(N-n)}(\phi_{\mu}(\theta_1,z_1))D_1^{2(N-m)}(D_1^{2m}(\Delta_{1,3})\Delta_{1,2})]\\& &+\sum_{m=0}^N\sum_{n=0}^{N-1}(-1)^{N+m+n}a_{\be,\gm,\lmd}^mb_{\lmd,\al,\mu}^nD_1^{2n+1}[D_1^{2(N-n)-1}(\phi_{\mu}(\theta_1,z_1))D_1^{2(N-m)}(D_1^{2m}(\Delta_{1,3})\Delta_{1,2})]\\& &+\sum_{m=0}^{N-1}\sum_{n=0}^N(-1)^{N+m+n}b_{\be,\gm,\lmd}^ma_{\lmd,\al,\mu}^nD_1^{2n}[D_1^{2(N-n)}(\phi_{\mu}(\theta_1,z_1))D_1^{2(N-m)-1}(D_1^{2m+1}(\Delta_{1,3})\Delta_{1,2})]\\& &+\sum_{m=0}^{N-1}\sum_{n=0}^{N-1}(-1)^{N+m+n+1}b_{\be,\gm,\lmd}^mb_{\lmd,\al,\mu}^nD_1^{2n+1}[D_1^{2(N-n)-1}(\phi_{\mu}(\theta_1,z_1))\\& &D_1^{2(N-m)-1}(D_1^{2m+1}(\Delta_{1,3})\Delta_{1,2})]\}\hspace{8.7cm}(5.19)\end{eqnarray*}
\begin{eqnarray*}&& [[\phi_{\gm}(\theta_3,z_3),\phi_{\al}(\theta_1,z_1)] ,\phi_{\be}(\theta_2,z_2)]\\&=&
z_2^{-1}z_3^{-1}\sum_{\lmd,\mu\in I}\{\sum_{m=0}^N\sum_{n=0}^N(-1)^ma_{\gm,\al,\lmd}^ma_{\lmd,\be,\mu}^nD_1^{2m}[D_1^{2(N-m)}[D_1^{2(N-n)}(\phi_{\mu}(\theta_1,z_1))D_1^{2n}(\Delta_{1,2})]\Delta_{1,3}]\\& &+\sum_{m=0}^N\sum_{n=0}^{N-1}(-1)^ma_{\gm,\al,\lmd}^mb_{\lmd,\be,\mu}^nD_1^{2m}[D_1^{2(N-m)}[D_1^{2(N-n)-1}(\phi_{\mu}(\theta_1,z_1))D_1^{2n+1}(\Delta_{1,2})]\Delta_{1,3}]\\& &+\sum_{m=0}^{N-1}\sum_{n=0}^N(-1)^{m+1}b_{\gm,\al,\lmd}^ma_{\lmd,\al,\mu}^nD_1^{2m+1}[D_1^{2(N-m)-1}[D_1^{2(N-n)}(\phi_{\mu}(\theta_1,z_1))D_1^{2n}(\Delta_{1,2})]\Delta_{1,3}]\\&&+\sum_{m=0}^{N-1}\sum_{n=0}^{N-1}(-1)^{m+1}b_{\gm,\al,\lmd}^mb_{\lmd,\al,\mu}^nD_1^{2m+1}[D_1^{2(N-m)-1}\\& &
[D_1^{2(N-n)-1}(\phi_{\mu}(\theta_1,z_1))D_1^{2n+1}(\Delta_{1,2})]\Delta_{1,3}]\}.\hspace{6.8cm}(5.20)\end{eqnarray*}

On the other hand, we have for $\bar{\xi}_1,\bar{\xi}_2,\bar{\xi}_3\in \Omega_0$: 
\begin{eqnarray*}& & \bar{\xi}_1(D_H\bar{\xi}_2H\bar{\xi}_3)\\&=&\sum_{\al,\lmd\in I}(D_H\bar{\xi}_2)_{\al,\lmd}(H\bar{\xi}_3)_{\lmd}\xi_{1\al}\\&=&\sum_{\al,\be,\gm,\lmd,\mu\in I}\{[\sum_{m=0}^Na_{\al,\be,\lmd}^mD^{2m}(\xi_{2\be})D^{2(N-m)}-\sum_{m=0}^{N-1}b_{\al,\be,\lmd}^mD^{2m+1}(\xi_{2\be})D^{2(N-m)-1}]\\&&[\sum_{n=0}^Na_{\lmd,\gm,\mu}^n\Phi_{\mu}(2(N-n)+1)D^{2n}(\xi_{3\gm})+\sum_{n=0}^{N-1}b_{\lmd,\gm,\mu}^n\Phi_{\mu}(2(N-n))D^{2n+1}(\xi_{3\gm})]\}\xi_{1\al}\\&=&\sum_{\al,\be,\gm,\lmd,\mu\in I}\{\sum_{m=0}^N\sum_{n=0}^Na_{\al,\be,\lmd}^ma_{\lmd,\gm,\mu}^nD^{2(N-m)}[\Phi_{\mu}(2(N-n)+1)D^{2n}(\xi_{3\gm})]D^{2m}(\xi_{2\be})\hspace{3cm}\end{eqnarray*}
\begin{eqnarray*}& &+\sum_{m=0}^N\sum_{n=0}^{N-1}a_{\al,\be,\lmd}^mb_{\lmd,\gm,\mu}^nD^{2(N-m)}[\Phi_{\mu}(2(N-n))D^{2n+1}(\xi_{3\gm})]D^{2m}(\xi_{2\be})\\& &-\sum_{m=0}^{N-1}\sum_{n=0}^Nb_{\al,\be,\lmd}^ma_{\lmd,\gm,\mu}^nD^{2(N-m)-1}[\Phi_{\mu}(2(N-n)+1)D^{2n}(\xi_{3\gm})]D^{2m+1}(\xi_{2\be})\\& &-\sum_{m=0}^{N-1}\sum_{n=0}^{N-1}b_{\al,\be,\lmd}^mb_{\lmd,\gm,\mu}^nD^{2(N-m)-1}[\Phi_{\mu}(2(N-n))D^{2n+1}(\xi_{3\gm})]D^{2m+1}(\xi_{2\be})\}\xi_{1\al},\hspace{1cm}(5.21)\end{eqnarray*}
\begin{eqnarray*}& & \bar{\xi}_2(D_H\bar{\xi}_3H\bar{\xi}_1)\\&=&\sum_{\al,\be,\gm,\lmd,\mu\in I}\{\sum_{m=0}^N\sum_{n=0}^N(-1)^{N+m+n}a_{\be,\gm,\lmd}^ma_{\lmd,\al,\mu}^nD^{2n}[\Phi_{\mu}(2(N-n)+1)D^{2(N-m)}[D^{2m}(\xi_{3\gm})\xi_{2\be}]]\\& &+\sum_{m=0}^N\sum_{n=0}^{N-1}(-1)^{N+m+n+1}a_{\be,\gm,\lmd}^mb_{\lmd,\al,\mu}^nD^{2n+1}[\Phi_{\mu}(2(N-n))D^{2(N-m)}[D^{2m}(\xi_{3\gm})\xi_{2\be}]]\\& &+\sum_{m=0}^{N-1}\sum_{n=0}^N(-1)^{N+m+n+1}b_{\be,\gm,\lmd}^ma_{\lmd,\al,\mu}^nD^{2n}[\Phi_{\mu}(2(N-n)+1)D^{2(N-m)-1}[D^{2m+1}(\xi_{3\gm})\xi_{2\be}]]\\& &+\sum_{m=0}^{N-1}\sum_{n=0}^{N-1}(-1)^{N+m+n}b_{\be,\gm,\lmd}^mb_{\lmd,\al,\mu}^nD^{2n+1}[\Phi_{\mu}(2(N-n))D^{2(N-m)-1}[D^{2m+1}(\xi_{3\gm})\xi_{2\be}]]\}\xi_{1\al},\hspace{0.5cm}(5.22)\end{eqnarray*}
\begin{eqnarray*}& & \bar{\xi}_3(D_H\bar{\xi}_1H\bar{\xi}_2)\\&=&
\sum_{\al,\be,\gm,\lmd,\mu\in I}\{\sum_{m=0}^N\sum_{n=0}^N(-1)^{m+1}a_{\gm,\al,\lmd}^ma_{\lmd,\be,\mu}^nD^{2m}[D^{2(N-m)}[\Phi_{\mu}(2(N-n)+1)D^{2n}(\xi_{2\be})]\xi_{3\gm}]\\& &+\sum_{m=0}^N\sum_{n=0}^{N-1}(-1)^{m+1}a_{\gm,\al,\lmd}^mb_{\lmd,\be,\mu}^nD^{2m}[D^{2(N-m)}[\Phi_{\mu}(2(N-n))D^{2n+1}(\xi_{2\be})]\xi_{3\gm}]\xi_{3\gm}]\\& &+\sum_{m=0}^{N-1}\sum_{n=0}^N(-1)^mb_{\gm,\al,\lmd}^ma_{\lmd,\be,\mu}^nD^{2m+1}[D^{2(N-m)-1}[\Phi_{\mu}(2(N-n)+1)D^{2n}(\xi_{2\be})]\\& &+\sum_{m=0}^{N-1}\sum_{n=0}^{N-1}(-1)^mb_{\gm,\al,\lmd}^mb_{\lmd,\be,\mu}^nD^{2m+1}[D^{2(N-m)-1}[\Phi_{\mu}(2(N-n))D^{2n+1}(\xi_{2\be})]\xi_{3\gm}]\}\xi_{1\al}
.\hspace{0.5cm}(5.23)\end{eqnarray*}
Comparing (5.18) and (5.21), (5.19) and (5.22), (5.20) and (5.23), we have:
\begin{eqnarray*}\hspace{1cm}& &[[\phi_{\al}(\theta_1,z_1),\phi_{\be}(\theta_2,z_2)] ,\phi_{\gm}(\theta_3,z_3)]+[[\phi_{\be}(\theta_2,z_2),\phi_{\gm}(\theta_3,z_3)] ,\phi_{\al}(\theta_1,z_1)]\\&=&-[[\phi_{\gm}(\theta_3,z_3),\phi_{\al}(\theta_1,z_1)] ,\phi_{\be}(\theta_2,z_2)].\hspace{6.9cm}(5.24)\end{eqnarray*}
Furthermore, we have:
\begin{eqnarray*}& &[[\phi_{\al}(\theta_1,z_1),\phi_{\be}(\theta_2,z_2)] ,\phi_{\gm}(\theta_3,z_3)]\\&=& [[\phi_{\al}^0(z_1),\phi_{\be}^0(z_2)] ,\phi_{\gm}^0(z_3)]\theta_1\theta_2\theta_3-[[\phi_{\al}^0(z_1),\phi_{\be}^1(z_2)] ,\phi_{\gm}^0(z_3)]\theta_1\theta_3\\&&+[[\phi_{\al}^1(z_1),\phi_{\be}^0(z_2)] ,\phi_{\gm}^0(z_3)]\theta_2\theta_3+[[\phi_{\al}^1(z_1),\phi_{\be}^1(z_2)] ,\phi_{\gm}^0(z_3)]\theta_3\hspace{6cm}\end{eqnarray*}
\begin{eqnarray*}&& +[[\phi_{\al}^0(z_1),\phi_{\be}^0(z_2)] ,\phi_{\gm}^1(z_3)]\theta_1\theta_2+[[\phi_{\al}^0(z_1),\phi_{\be}^1(z_2)] ,\phi_{\gm}^1(z_3)]\theta_1\\&&-[[\phi_{\al}^1(z_1),\phi_{\be}^0(z_2)] ,\phi_{\gm}^1(z_3)]\theta_2+[[\phi_{\al}^1(z_1),\phi_{\be}^1(z_2)] ,\phi_{\gm}^1(z_3)],\hspace{4.5cm}(5.25)\end{eqnarray*}
\begin{eqnarray*}& &[[\phi_{\be}(\theta_2,z_2),\phi_{\gm}(\theta_3,z_3)] ,\phi_{\al}(\theta_1,z_1)]\\&=&  [[\phi_{\be}^0(z_2),\phi_{\gm}^0(z_3)] ,\phi_{\al}^0(z_1)]\theta_1\theta_2\theta_3-[[\phi_{\be}^1(z_2),\phi_{\gm}^0(z_3)] ,\phi_{\al}^0(z_1)]\theta_1\theta_3\\&&+[[\phi_{\be}^0(z_2),\phi_{\gm}^0(z_3)] ,\phi_{\al}^1(z_1)]\theta_2\theta_3-[[\phi_{\be}^1(z_2),\phi_{\gm}^0(z_3)] ,\phi_{\al}^1(z_1)]\theta_3\\& &+[[\phi_{\be}^0(z_2),\phi_{\gm}^1(z_3)] ,\phi_{\al}^0(z_1)]\theta_1\theta_2+[[\phi_{\be}^1(z_2),\phi_{\gm}^1(z_3)] ,\phi_{\al}^0(z_1)]\theta_1\\& &+[[\phi_{\be}^0(z_2),\phi_{\gm}^1(z_3)] ,\phi_{\al}^1(z_1)]\theta_2+[[\phi_{\be}^1(z_2),\phi_{\gm}^1(z_3)] ,\phi_{\al}^1(z_1)],\hspace{4.1cm}(5.26)\end{eqnarray*}
\begin{eqnarray*}& &[[\phi_{\gm}(\theta_3,z_3),\phi_{\al}(\theta_1,z_1)] ,\phi_{\be}(\theta_2,z_2)]\\&=& [[\phi_{\gm}^0(z_3),\phi_{\al}^0(z_1)] ,\phi_{\be}^0(z_2)]\theta_1\theta_2\theta_3- [[\phi_{\gm}^0(z_3),\phi_{\al}^0(z_1)] ,\phi_{\be}^1(z_2)]\theta_1\theta_3\\
&& + [[\phi_{\gm}^0(z_3),\phi_{\al}^1(z_1)] ,\phi_{\be}^0(z_2)]\theta_2\theta_3+ [[\phi_{\gm}^0(z_3),\phi_{\al}^1(z_1)] ,\phi_{\be}^1(z_2)]\theta_3\\&&+ [[\phi_{\gm}^1(z_3),\phi_{\al}^0(z_1)] ,\phi_{\be}^0(z_2)]\theta_1\theta_2-[[\phi_{\gm}^1(z_3),\phi_{\al}^0(z_1)] ,\phi_{\be}^1(z_2)]\theta_1\\& &+ [[\phi_{\gm}^1(z_3),\phi_{\al}^1(z_1)] ,\phi_{\be}^0(z_2)]\theta_2+ [[\phi_{\gm}^1(z_3),\phi_{\al}^1(z_1)] ,\phi_{\be}^1(z_2)].\hspace{4.1cm}(5.27)\end{eqnarray*}
Now (5.24-27) imply:
\begin{eqnarray*}& &[[\phi_{\al}^i(z_1),\phi_{\be}^j(z_2)],\phi_{\gm}^k(z_3)]+(-1)^{i(j+k)}[[\phi_{\be}^j(z_2),\phi_{\gm}^k(z_3)],\phi_{\al}^i(z_1)]\\&=&-(-1)^{k(i+j)}[[\phi_{\al}^k(z_3),\phi_{\al}^i(z_1)],\phi_{\be}^j(z_2)]\hspace{7.6cm}(5.28)\end{eqnarray*}
which is the super Jacobi identity. $\qquad\Box$
\psp

{\bf Example}. If we let $a_{\al,\be}=0$ in (3.1), then the Hamiltonian operator (3.1) is a special case of (5.1). Therefore, each NX-bialgebra induces a Lie superalgebra. It can be proved that $\{a_{\al,\be}\mid \al,\be \in I\}$ induces
a one-dimeional central extension of the Lie superalgebra. In general, Threorem 5.3 holds for any linear Hamiltonian superoperator of type 1. 

  Let $0<n\in \Bbb{Z}$ and let $L$ be a vector space with a basis $\{\phi_i(n),c\mid n\in \Bbb{Z}/2, i=0,1,...,n-1\}$. Besides using (5.8-10) and (5.14), we also assume:
$$\theta_ic=c\theta_i,\qquad cz_i=z_ic\qquad\qquad\mbox{for}\;\;i=1,2.\eqno(5.29)$$
By the example in Section 3, we have the following type-1 Hamiltonian superoperator $H$:
$$H^1_{i,j}=H_{i,j}^0=\delta_{i,0}\delta_{j,0}D^5+(i+j+3)\Phi_{i+j}D^2+\Phi_{i+j}(2)D+(j+2)\Phi_{i+j}(3)\eqno(5.30)$$
for $i,j=0,1,...,n-1$. Here we have used the convention that $\Phi_l=0$ if $l\geq n$. We define the opration $[\cdot,\cdot]$ on $L$ by:
$$[u,c]=[c,u]=0\qquad\qquad\mbox{for}\;\;u\in L\eqno(5.31)$$
and
\begin{eqnarray*}&& [\phi_i(\theta_1,z_1),\phi_j(\theta_2,z_2)]\\&=& z_2^{-1}\{
\delta_{i,0}\delta_{j,0}D^5\Delta_{1,2}c+(i+j+3)\Phi_{i+j}(\theta_1,z_1)D_1^2\Delta_{1,2}\\& &+D_1(\phi_{i+j}(\theta_1,z_1))D_1\Delta_{1,2}+(j+2)D_1^2(\phi_{i+j}(\theta_1,z_1))\Delta_{1,2}\}\\&=& z_2^{-1}\{\delta_{i,0}\delta_{j,0}[\partial_{z_1}^2\delta(z_1/z_2)c-\partial_{z_1}^3\delta(z_1/z_2)c\theta_1\theta_2]+(i+j+3)[\phi_{i+j}^0(z_1)\theta_1+\phi_{i+j}^1(z_1)]\\& &
(\theta_1-\theta_2)\partial_{z_1}\delta(z_1/z_2)+[\phi_{i+j}^0(z_1)-\partial_{z_1}(\phi_{i+j}^1(z_1))\theta_1][\delta(z_1/z_2)-\partial_{z_1}\delta(z_1/z_2)\theta_1\theta_2]\\& &+(j+2)[\partial_{z_1}(\phi_{i+j}^0(z_1))\theta_1+\partial_{z_1}(\phi_{i+j}^1(z_1))](\theta_1-\theta_2)\delta(z_1/z_2)\}\\&=&z_2^{-1}\{\delta_{i,0}\delta_{j,0}\partial_{z_1}^2\delta(z_1/z_2)c+\phi_{i+j}^0(z_1)\delta(z_1/z_2)\\&&+[(i+j+3)\phi_{i+j}^1(z_1)\partial_{z_1}\delta(z_1/z_2)+(j+1)\partial_{z_1}(\phi_{i+j}^1(z_1))\delta(z_1/z_2)]\theta_1\\& &-[(i+j+3)\phi_{i+j}^1(z_1)\partial_{z_1}\delta(z_1/z_2)+(j+2)\partial_{z_1}(\phi_{i+j}^1(z_1))\delta(z_1/z_2)]\theta_2\\& &-[\delta_{i,0}\delta_{j,0}\partial_{z_1}^3\delta(z_1/z_2)c+(i+j+4)\phi_{i+j}^0(z_1)\partial_{z_1}\delta(z_1/z_2)\\& &+(j+2)\partial_{z_1}(\phi_{i+j}^0(z_1))\delta(z_1/z_2)]\theta_1\theta_2\}\hspace{7.6cm}(5.32)\end{eqnarray*}
for $i,j=0,1,...,n-1.$ Thus we have
$$[\phi_i^1(z_1),\phi_j^1(z_2)]=z_2^{-1}[\delta_{i,0}\delta_{j,0}\partial_{z_1}^2\delta(z_1/z_2)c+\phi_{i+j}^0(z_1)\delta(z_1/z_2)],\eqno(5.33)$$
$$-[\phi_i^0(z_1),\phi_j^1(z_2)]=z_2^{-1}[(i+j+3)\phi_{i+j}^1(z_1)\partial_{z_1}\delta(z_1/z_2)+(j+1)\partial_{z_1}(\phi_{i+j}^1(z_1))\delta(z_1/z_2)],\eqno(5.34)$$
$$[\phi_i^1(z_1),\phi_j^0(z_2)]=-z_2^{-1}[(i+j+3)\phi_{i+j}^1(z_1)\partial_{z_1}\delta(z_1/z_2)+(j+2)\partial_{z_1}(\phi_{i+j}^1(z_1))\delta(z_1/z_2)],\eqno(5.35)$$
\begin{eqnarray*}& &[\phi_i^0(z_1),\phi_j^0(z_2)]\\&=&-z_2^{-1}[\delta_{i,0}\delta_{j,0}\partial_{z_1}^3\delta(z_1/z_2)c+(i+j+4)\phi_{i+j}^0(z_1)\partial_{z_1}\delta(z_1/z_2)\\& &+(j+2)\partial_{z_1}(\phi_{i+j}^0(z_1))\delta(z_1/z_2)].\hspace{8.4cm}(5.36)\end{eqnarray*}
Note that (5.33-36) are equivalent to:
$$\left[\phi_i\left(m+{1\over 2}\right),\phi_j\left(n+{1\over 2}\right)\right]=\delta_{i,0}\delta_{j,0}\delta_{m+n+1,0}(n+1)nc+\phi_{i+j}(m+n+1),\eqno(5.37)$$
$$\left[\phi_i\left(m+{1\over 2}\right),\phi_j(n)\right]=[(j+2)(m+1)-(i+1)(n+1)]\phi_{i+j}\left(m+n+{1\over 2}\right),\eqno(5.38)$$
\begin{eqnarray*}[\phi_i(m),\phi_j(n)]&=&-\delta_{i,0}\delta_{j,0}\delta_{m+n+1,0}(n+1)n(n-1)c\\& &+[(j+2)(m+1)-(i+2)(n+1)]\phi_{i+j}(m+n)\hspace{3cm}(5.39)\end{eqnarray*}
for $i,j=0,1,...,n-1;\;m,n\in \Bbb{Z}$. Therefore, we obtain a Lie superalgebra
$(L,[\cdot,\cdot])$ that is a natural generalization of the Super-Virasoro algebra.
\psp

{\bf Remark 5.4}. (a) Lie superalgebras induced by Novikov-Poisson algebras whose Novikov algebras are simple were studied in [X5].

(b) We still do not know how to connect linear Hamiltonian superoperators of type 0 with Lie superalgebras.

\hspace{0.5cm}

\noindent{\Large \bf References}

\hspace{0.5cm}

\begin{description}

\item[{[BN]}] A. A. Balinskii and S. P. Novikov, Poisson brackets of hydrodynamic type, Frobenius algebras and Lie algebras, {\it Soviet Math. Dokl.} Vol. {\bf 32} (1985), No. {\bf 1}, 228-231.

\item[{[DGM]}] L. Dolan, P. Goddard and P. Montague, Conformal field theory of twisted vertex operators, {\it Nucl. Phys.} {\bf B338} (1990) 529-601.

\item[{[Da1]}]
Yu. L. Daletsky, Lie superalgebras in Hamiltonian operator theory, In: {\it Nonlinear and Turbulent Processes in Physics}, ed. V. E. Zakharov, 19984, pp. 1307-1312.

\item[{[Da2]}]
Yu. L. Daletsky, Hamiltonian operators in graded formal calculus of variables,
{\it Func. Anal. Appl.} {\bf 20} (1986), 136-138.

\item[{[De]}] B. DeWitt, {\it Supermanifolds,} Second Edition, Cambridge University Press, 1992. 

\item[{[FFR]}]
A. J. Feingold, I. B. Frenkel and J. F. Ries, Spinor construction of vertex operator algebras, triality and $E_8^{(1)}$, {\it Contemp. Math.} {\bf 121}, 1991.

\item[{[FLM]}] 
I. B. Frenkel, J. Lepowsky and A. Meurman, {\it Vertex Operator
Algebras and the Monster}, Pure and Applied Math. Academic Press, 1988.

\item[{[GDi1]}] I. M. Gel'fand and L. A. Dikii, Asymptotic behaviour of the resolvent of Sturm-Liouville equations and the algebra of the Korteweg-de Vries equations, {\it Russian Math. Surveys} {\bf 30:5} (1975), 77-113.

\item[{[GDi2]}] 
I. M. Gel'fand and L. A. Dikii, A Lie algebra structure in a formal variational Calculation, {\it Func. Anal. Appl.}  {\bf 10} (1976), 16-22.

\item[{[GDo]}] 
I. M. Gel'fand and I. Ya. Dorfman, Hamiltonian operators and algebraic structures related to them, {\it Func. Anal. Appl.}  {\bf 13} (1979),  248-262.

\item[{[M]}] P. Mathieu, Supersymetry extension of the Korteweg-de Vries equation, {\it J. Math. Phys.} {\bf 29}(11) (1988), 2499-2507.

\item[{[NO]}] J. W. Negele and H. Orland, {\it Quantum many-particle systems}, Addison-Wesley Publishing Company, 1988.

\item[{[O1]}]
J. Marshall Osborn, Novikov algebras, {\it Nova J. Algebra} \& {\it Geom.} {\bf 1} (1992), 1-14.

\item[{[O2]}]
J. Marshall Osborn, Simple Novikov algebras with an idempotent, {\it Comm. Algebra} {\bf 20} (1992), No. 9, 2729-2753.

\item[{[O3]}]
J. Marshall Osborn, Infinite dimensional Novikov algebras of characteristic 0, {\it J. Algebra} {\bf 167} (1994), 146-167.

\item[{[O4]}]
J. Marshall Osborn, Modules for Novikov algebras, {\it Proceeding of the II International Congress on Algebra, Barnaul, 1991.}

\item[{[O5]}]
J. Marshall Osborn, Modules for Novikov algebras of characteristic 0, {\it preprint}.

\item[{[T]}] 
H. Tsukada,  Vertex operator superalgebras, {\it Comm. Math. Phys.} {\bf 18} (1990), 2249-2274.

\item[{[X1]}]
X. Xu, On spinor vertex operator algebras and their modules, {\it J. Algebra} {\bf 191}, 427-460.

\item[{[X2]}], X. Xu, Hamiltonian operators and associative algebras with a derivation, {\it  Lett. Math. Phys.} {\bf 33} (1995), 1-6.

\item[{[X3]}] X. Xu, Hamiltonian superoperators, {\it J. Phys A: Math. Gen.} {\bf 28} (1995), 1681-1698.

\item[{[X4]}] X. Xu, On simple Novikov algebras and their irreducible modules, {\it J. Algebra} {\bf 185} (1996), 905-934.

\item[{[X5]}] X. Xu, Novikov-Poisson Algebras, {\it J. Algebra} {\bf 190} (1997), 253-279.

\item[{[X6]}] X. Xu, Skew-symmetric differential operators and combinatorial identities, {\it Mh. Math} {\bf 127} (1999), 243-258.

\item[{[Z]}]
E. I. Zel'manov, On a class of local translation invariant Lie algebras, {\it Soviet Math. Dokl.} Vol {\bf 35} (1987), No. {\bf 1}, 216-218.

\end{description}
\end{document}